\documentclass[11pt]{article}%
\usepackage{graphicx}
\usepackage{latexsym}
\usepackage{amssymb}
\usepackage{amsmath}
\usepackage{layout}
\usepackage{amsfonts}%
\newtheorem{proposition}{Proposition}
\newtheorem{lemma}{Lemma}
\newtheorem{definition}{Definition}

\newtheorem{theorem}{Theorem}

\newenvironment{proof}[1][Proof]{\textbf{#1.} }{\ \rule{0.5em}{0.5em}}
\ifx\pdfoutput\relax\let\pdfoutput=\undefined\fi
\newcount\msipdfoutput
\ifx\pdfoutput\undefined\else
\ifcase\pdfoutput\else
\ifx\paperwidth\undefined\else
\ifdim\paperheight=0pt\relax\else\pdfpageheight\paperheight\fi
\ifdim\paperwidth=0pt\relax\else\pdfpagewidth\paperwidth\fi
\fi\fi\fi
\begin{document}

\title{Variations and Hurst index estimation for a Rosenblatt process using longer filters}
\author{Alexandra Chronopoulou$^{1,}\thanks{Authors partially supported by NSF grant
0606615}\qquad$ Ciprian A. Tudor$^{2,} \thanks{Associate member: SAMOS-MATISSE, Centre d'Economie de La Sorbonne, Universit\'{e} de Paris 1 Panth\'eon-Sorbonne ,90, rue de Tolbiac, 75634, Paris, France.} \qquad$Frederi G. Viens$^{1,*
,}\thanks{Corresponding author}$\vspace*{0.1in}\\$^{1}$ Department of Statistics, Purdue University, \\150 N. University St., West Lafayette, IN 47907-2067, USA. \\achronop@purdue.edu\quad\quad viens@purdue.edu $^{\dag}$\vspace*{0.1in} \\$^{2}$Laboratoire Paul Painlev\'e, Universit\'e de Lille 1,\\F-59655 Villeneuve d'Ascq, France. \\tudor@math.univ-lille1.fr\vspace*{0.1in}}
\maketitle

\begin{abstract}
The Rosenblatt process is a self-similar non-Gaussian process which lives in
second Wiener chaos, and occurs as the limit of correlated random sequences in
so-called \textquotedblleft non-central limit theorems\textquotedblright. It
shares the same covariance as fractional Brownian motion. We study the
asymptotic distribution of the quadratic variations of the Rosenblatt process
based on long filters, including filters based on high-order finite-difference
and wavelet-based schemes. We find exact formulas for the limiting
distributions, which we then use to devise strongly consistent estimators of
the self-similarity parameter $H$. Unlike the case of fractional Brownian
motion, no matter now high the filter orders are, the estimators are never
asymptotically normal, converging instead in the mean square to the observed
value of the Rosenblatt process at time $1$.

\end{abstract}

\vskip0.3cm

\textbf{2000 AMS Classification: }Primary: 60G18; Secondary 60F05, 60H05, 62F12.

\vskip0.2cm


\section{Introduction}

Self-similar stochastic processes are of practical interest in various
applications, including econometrics, internet traffic, and hydrology. These
are processes $X=\left\{  X\left(  t\right)  :t\geq0\right\}  $ whose
dependence on the time parameter $t$ is self-similar, in the sense that there
exists a (self-similarity) parameter $H\in(0,1)$ such that for any constant
$c\geq0$, $\left\{  X\left(  ct\right)  :t\geq0\right\}  $ and $\left\{
c^{H}X\left(  t\right)  :t\geq0\right\}  $ have the same distribution. These
processes are often endowed with other distinctive properties.

The fractional Brownian motion (fBm) is the usual candidate to model phenomena
in which the selfsimilarity property can be observed from the empirical data.
This fBm $B^{H}$ is the continuous centered Gaussian process with covariance
function
\begin{equation}
R^{H}(t,s):=\mathbf{E}\left[  B^{H}\left(  t\right)  B^{H}\left(  s\right)
\right]  =\frac{1}{2}(t^{2H}+s^{2H}-|t-s|^{2H}). \label{cov}%
\end{equation}
The parameter $H\ $characterizes all the important properties of the process.
In addition to being self-similar with parameter $H$, which is evident from
the covariance function, fBm has correlated increments: in fact, from
(\ref{cov}) we get, as $n\rightarrow\infty$,
\begin{equation}
\mathbf{E}\left[  \left(  B^{H}\left(  n\right)  -B^{H}\left(  1\right)
\right)  B^{H}\left(  1\right)  \right]  =H\left(  2H-1\right)  n^{2H-2}%
+o\left(  n^{2H-2}\right)  ; \label{long}%
\end{equation}
when $H<1/2$, the increments are negatively correlated and the correlation
decays more slowly than quadratically; when $H>1/2$, the increments are
positively correlated and the correlation decays so slowly that they are not
summable, a situation which is commonly known as the long memory property. The
covariance structure (\ref{cov}) also implies%
\begin{equation}
\mathbf{E}\left[  \left(  B^{H}\left(  t\right)  -B^{H}\left(  s\right)
\right)  ^{2}\right]  =\left\vert t-s\right\vert ^{2H}; \label{canon}%
\end{equation}
this property shows that the increments of fBm are stationary and
self-similar; its immediate consequence for higher moments can be used, via
the so-called Kolmogorov continuity criterion, to imply that $B^{H}$ has paths
which are almost-surely ($H-\varepsilon$)-H\"{o}lder-continuous for any
$\varepsilon>0$.

It turns out that fBm is the \emph{only} continuous Gaussian process which is selfsimilar with
stationary increments. However, there are many more stochastic processes which, except for the
Gaussian character, share all the other properties above for $H>1/2$ (i.e. (\ref{cov}) which
implies (\ref{long}), the long-memory property, (\ref{canon}), and in many cases the
H\"{o}lder-continuity). In some models the Gaussian assumption may be implausible and in this
case one needs to use a different self-similar process with stationary increments to model the
phenomenon. Natural candidates are the Hermite processes: these non-Gaussian stochastic
processes appear as limits in the so-called Non-Central Limit Theorem (see \cite{BrMa},
\cite{DM}, \cite{Ta1}) and do indeed have all the properties listed above. While fBm can be
expressed as a Wiener integral with respect to the standard Wiener process, i.e. the integral of
a deterministic kernel w.r.t. a standard Brownian motion, the Hermite process of order $q\geq2$
is a $q$th iterated integral of a deterministic function with $q$ variables with respect to a
standard Brownian motion. When $q=2$, the Hermite process is called the Rosenblatt process. This
stochastic process typically appears as a limiting model in various applications such as unit
the root testing problem (see \cite{Wu} ), semiparametric approach to hypothesis test (see
\cite{Hall}), or long-range dependence estimation (see \cite{KG}). On the other hand, since it
is non-Gaussian and self-similar with stationary increments, the Rosenblatt process can also be
an input in models where self-similarity is observed in empirical data which appears to be
non-Gaussian. The need of non-Gaussian self-similar processes in practice (for example in
hydrology) is mentioned in the paper \cite{Taqqu3} based on the study of stochastic modeling for
river-flow time series in \cite{LK}.  Recent interest in the Rosenblatt and other Hermite
processes, due in part to their non-Gaussian character, and in part for their independent
mathematical value, is evidenced by the following references: \cite{BN}, \cite{CNT},
\cite{Es-Tu}, \cite{MaTu}, \cite{NNT}, \cite{MPR}, \cite{T}, \cite{TV}.

The results in these articles, and in the previous references on the
non-central limit theorem, have one point in common: of all the Hermite
processes, the most important one in terms of limit theorem, apart from fBm,
is the Rosenblatt process. As such, it should be the first non-Gaussian
self-similar process for which to develop a full statistical estimation
theory. This is one motivation for writing this article.

Since the Hurst parameter $H$, thus called in reference to the hydrologist who
discovered its original practical importance (see \cite{Hurst}), characterizes
all the important properties of a Hermite process, its proper statistical
estimation is of the utmost importance. Several statistics have been
introduced to this end in the case of fBm, such as variograms, maximum
likelihood estimators, or spectral methods, k-variations and wavelets.
Information on these various approaches, apart from wavelets, for fBm and
other long-memory processes, can be found in the book of Beran \cite{Beran}.
More details about the wavelet-based approach can be found in \cite{BaTu},
\cite{Fl} and \cite{WO}.

In this article, we will concentrate on one of the more popular methods to
estimate $H$: the study of power variations; it is particularly well-adapted
to the non-Gaussian Hermite processes, because explicit calculations can be
performed via Wiener chaos analysis. In its simplest form, the $k$th power
variation statistic of a process $\left\{  X_{t}:t\in\lbrack0,1]\right\}  $,
calculated using $N$ data points, is defined as following quantity (the
absolute value of the increment may be used in the definition for non-even
powers):%
\begin{equation}
V_{N}:=\frac{1}{N}\left[  \sum_{i=0}^{N-1}\frac{\left(  X_{\frac{i+1}{N}%
}-X_{\frac{i}{N}}\right)  ^{k}}{\mathbf{E}\left(  X_{\frac{i+1}{N}}%
-X_{\frac{i}{N}}\right)  ^{k}}-1\right]  . \label{vp}%
\end{equation}

There exists a direct connection between the behavior of the variations and
the convergence of an estimator for the selfsimilarity order based on these
variations (see \cite{coeur}, \cite{TV}): if the renormalized variation
satisfies a central limit theorem then so does the estimator, a desirable fact
for statistical purposes.

The recent paper \cite{TV} studies the quadratic variation of the Rosenblatt
process $Z$ (the $V_{N}$ above with $k=2$), exhibiting the following facts:
the normalized sequence $N^{1-H}V_{N}$ satisfies a non-central limit theorem,
it converges in the mean square to the Rosenblatt random variable $Z\left(
1\right)  $ (value of the process $Z$ at time $1$); from this, we can
construct an estimator for $H$ whose behavior is still non-normal. The same
result is also obtained in the case of the estimators based on the wavelet
coefficients (see \cite{BaTu}). In the simpler case of fBm, this situation
still occurs when $H>3/4$ (see for instance \cite{TV2}). For statistical
applications, a situation in which asymptotic normality holds might be
preferable. To achieve this for fBm, it has been known for some time that one
may use \textquotedblleft longer filters\textquotedblright\ (that means,
replacing the increments $X_{\frac{i+1}{N}}-X_{\frac{i}{N}}$ by the second-order increments $X_{\frac
{i+1}{N}}-2X_{\frac{i}{N}}+X_{\frac{i-1}{N}}$, or higher order increments for
instance; see \cite{coeur}). To have asymptotic normality in the case of the
Rosenblatt process, it was shown in \cite{TV} that one may perform a
compensation of the non-normal component of the quadratic variation. In fact,
this is possible only in the case of the Rosenblatt process; it is not
possible for higher-order Hermite processes, and is not possible for fBm with
$H>3/4$ [recall that the case of fBm with $H\leq3/4$ does not require any
compensation]. The compensation technique for the Rosenblatt process yields
asymptotic variances which are difficult to calculate and may be very high.

The question then arises to find out whether using longer filters for the
Rosenblatt process might yield asymptotically normal estimators, and/or might
result in low asymptotic variances. In this article, using recent results on
limit theorems for multiple stochastic integrals based on the Malliavin
calculus (see \cite{NP}, \cite{NOT}), we will see that the answer to the first
question is negative, while the answer to the second question is affirmative.
We will use quadratic variations ($k=2$) for simplicity. A summary of our
results is as follows. Here $\Omega$ denotes the underlying probability space,
and $L^{1}\left(  \Omega\right)  $ and $L^{2}\left(  \Omega\right)  $ are the
usual spaces of integrable and square-integrable random variables.

\begin{itemize}
\item $V_{N}=T_{2}+T_{4}$ where $T_{i}$ is in the $i$th Wiener chaos
(Proposition \ref{VNchaos}).

\item $\frac{\sqrt{N}}{c_{1,H}}T_{4}$ converges in distribution to a standard
normal (Theorem \ref{thmT4}), where $c_{1,H}$ is given in Proposition
\ref{propET4}.

\item $\frac{N^{1-H}}{\sqrt{c_{2,H}}}V_{N}$ and $\frac{N^{1-H}}{\sqrt{c_{2,H}%
}}T_{2}$ both converge in $L^{2}(\Omega)$ to the Rosenblatt random variable
$Z(1)$ (Theorem \ref{thmV}); the asymptotic variance $c_{2,H}$ is given
explicitly in formula (\ref{c2}) in Proposition \ref{propc2}.

\item There exists a strongly consistent estimator $\hat{H}_{N}$ for $H$ based
on $V_{N}$ (Theorem \ref{thmconsistent}), and $2\;c_{2,H}^{-1/2}(\log
N)N^{1-\hat{H}_{N}}\left(  \hat{H}_{N}-H\right)  $ converges in $L^{1}\left(
\Omega\right)  $ to a Rosenblatt random variable (Theorem \ref{ThmStat}). Here
$c_{2,H}$ is again given in (\ref{c2}). Note that while the rate of
convergence of the estimator, of order $N^{-1+H}\log^{-1}N$, depends on $H$,
the convergence result above can be used without knowledge of $H$ since one
may plug in $\hat{H}_{N}$ instead of $H$ in the convergence rate.

\item The asymptotic variance $c_{2,H}$ in the above convergence decreases as
the length of the filter increases; this decrease is much faster for
wavelets-based filters than for finite-difference-based filters: for values of
$H<0.95$, $c_{2,H}$ reaches values below $5\%$ for wavelet filters of length
less than 6, but for finite-difference filters of length no less than 16.

\item When $H\in(1/2,2/3)$, then $\frac{N}{c_{3,H}}\left[  V_{N}-\frac
{\sqrt{c_{2,H}}}{N^{1-H}}Z(1)\right]  $ converges in distribution to a
standard normal, where $c_{2,H}$ is given explicitly in formula (\ref{c2}) and
$c_{3,H}$ in formula (\ref{c3}). Similarly, for the estimator we have that
$\frac{N}{c_{3,H}}\left[  -2\log(\hat{H}_{N}-H)-\frac{\sqrt{c_{2,H}}}{N^{1-H}%
}Z(1)\right]  $ converges in distribution to the same standard normal.
However, no mater how much we increase the order and/or the length of the
filter, we cannot improve the threshold of 2/3 for $H$.
\end{itemize}

What prevents the normalization of $V_{N}$ from converging to a Gaussian, no
matter how long the filter is, is the distinction between the two terms
$T_{2}$ and $T_{4}$. In the case of fractional Brownian motion, $V_{N}$
contains only one \textquotedblleft$T_{2}$\textquotedblright-type term (second
chaos), but this term has a behavior similar to our term $T_{4}$, and does
converge to a normal when the filter is long enough; this fact has been noted
before (see \cite{coeur}). In our case, the normalized $T_{2}$ always
converges (in $L^{2}\left(  \Omega\right)  $) to a Rosenblatt random variable;
the piece that sometimes has normal asymptotics is $T_{4}$, but since $T_{2}$
always dominates it, $V_{N}$'s behavior is always that of $T_{2}$. This sort
of phenomenon was already noted in \cite{CNT} with the order-one filter for
all non-Gaussian Hermite processes, but now we know it occurs also for the
simplest Hermite process that is not fBm, for filters of all orders.\vspace
*{0.1in}

The organization of our paper is as follows. Section 2 summarizes the
stochastic analytic tools we will use, and gives the definitions of the
Rosenblatt process and the filter variations. Therein we also establish a
specific representation of the 2-power variation as the sum of two terms, one
in the second Wiener chaos, which we call $T_{2}$, and another, $T_{4}$, in
the fourth Wiener chaos. Section 3 establishes the correct normalizing factors
for the variations, by computing second moments, showing in particular that
$T_{2}$ is the dominant term. Section 4 proves that the renormalized $T_{4}$
is asymptotically normal. Section 5 proves that $T_{2}$ converges in
$L^{2}\left(  \Omega\right)  $ to the value $Z\left(  1\right)  $ of the
Rosenblatt process at time $1$. In Section 6 it is shown that the variation
obtained by subtracting this observed limit of $T_{2}$ leads to a correction
term which is asymptotically normal. Section 7 establishes the strong
consistency of the estimator $\hat{H}$ for $H$ based on the variations, and
proves that the renormalized estimator converges to a Rosenblatt random
variable in $L^{1}\left(  \Omega\right)  $. Its asymptotic variance is given
explicitly for any filter, thanks to the calculations in Section 3. In Section
8, we compare the numerical values of the asymptotic variances for various
choices of filters, including finite-difference filters and wavelet-based
filters, concluding that the latter are more efficient.

\section{Preliminaries}

\label{Prem}

\subsection{Basic tools on multiple Wiener-It\^{o} integrals}

Let $\left\{  W_{t}:t\in\right[  0,1]\}$ be a classical Wiener process on a
standard Wiener space $\left(  \Omega,{\mathcal{F}},P\right)  $. If a
symmetric function $f\in L^{2}([0,1]^{n})$ is given, the multiple
Wiener-It\^{o} integral $I_{n}\left(  f\right)  $ of $f$ with respect to $W$
is constructed and studied in detail in \cite[Chapter 1]{N}. Here we collect
the results we will need. For the most part, the results in this subsection
will be used in the technical portions of our proofs, which are in the
Appendix. One can construct the multiple integral starting from simple
functions of the form $f:=\sum_{i_{1},\ldots,i_{n}}c_{i_{1},\ldots i_{n}%
}1_{A_{i_{1}}\times\ldots\times A_{i_{n}}}$ where the coefficient $c_{i_{1},..,i_{n}}$ is zero if two indices are equal and the sets $A_{i_{j}}$ are
disjoint intervals, by setting%
\[
I_{n}(f):=\sum_{i_{1},\ldots,i_{n}}c_{i_{1},\ldots i_{n}}W(A_{i_{1}})\ldots
W(A_{i_{n}})
\]
where we put $W\left(  1_{[a,b]}\right)  =W([a,b])=W_{b}-W_{a}$; then the
integral is extended to all symmetric functions in $L^{2}([0,1]^{n})$ by a
density argument. It is also convenient to note that this construction
coincides with the iterated It\^{o} stochastic integral%
\[
I_{n}(f)=n!\int_{0}^{1}\int_{0}^{t_{n}}\ldots\int_{0}^{t_{2}}f(t_{1}%
,\ldots,t_{n})dW_{t_{1}}\ldots dW_{t_{n}}.
\]
The application $I_{n}$ is extended to non-symmetric functions $f$ via%
\begin{equation}
I_{n}(f)=I_{n}\big(\tilde{f}\big) \label{ftilde}%
\end{equation}
where $\tilde{f}$ denotes the symmetrization of $f$ defined by $\tilde
{f}(x_{1},\ldots,x_{x})=\frac{1}{n!}\sum_{\sigma\in\mathcal{S}_{n}}%
f(x_{\sigma(1)},\ldots,x_{\sigma(n)})$. The map $\left(  n!\right)
^{-1/2}I_{n}$ can then be seen to be an isometry from $L^{2}([0,1]^{n})$ to
$L^{2}(\Omega)$. The $n$th Wiener chaos is the set of all integrals $\left\{
I_{n}\left(  f\right)  :f\in L^{2}([0,1]^{n})\right\}  $; the Wiener chaoses
form orthogonal sets in $L^{2}\left(  \Omega\right)  $. Summarizing, we have%
\begin{align}
\mathbf{E}\left(  I_{n}(f)I_{m}(g)\right)   &  =n!\langle f,g\rangle
_{L^{2}([0,1]^{n})}\quad\mbox{if }m=n,\label{isom}\\
\mathbf{E}\left(  I_{n}(f)I_{m}(g)\right)   &  =0\quad\mbox{if }m\not =%
n.\nonumber
\end{align}

The product for two multiple integrals can be expanded explicitly (see
\cite{N}): if $f\in L^{2}([0,1]^{n})$ and $g\in L^{2}([0,1]^{m})$ are
symmetric, then it holds that
\begin{equation}
I_{n}(f)I_{m}(g)=\sum_{\ell=0}^{m\wedge n}\ell!C_{m}^{\ell}C_{n}^{\ell
}I_{m+n-2\ell}(f\otimes_{\ell}g) \label{prod}%
\end{equation}
where the contraction $f\otimes_{\ell}g$ belongs to $L^{2}([0,1]^{m+n-2\ell})$
for $\ell=0,1,\ldots,m\wedge n$ and is given by
\begin{eqnarray*}
&& (f\otimes_{\ell}g)(s_{1},\ldots,s_{n-\ell},t_{1},\ldots,t_{m-\ell})\\
&=& \int_{[0,1]^{\ell}}f(s_{1},\ldots,s_{n-\ell},u_{1},\ldots,u_{\ell}
)g(t_{1},\ldots,t_{m-\ell},u_{1},\ldots,u_{\ell})du_{1}\ldots du_{\ell}.
\end{eqnarray*}
Note that the contraction $(f\otimes_{\ell} g) $ is not necessary symmetric. We will denote by $(f\tilde{\otimes }_{\ell} g)$ its symmetrization.
\bigskip

Our analysis will be based on the following result, due to Nualart and
Peccati (see Theorem 1 in \cite{NP}).

\begin{proposition}
\label{t2}Let $n$ be a fixed integer. Let $I_{n}(f_{N})$ be a sequence
of symmetric square integrable random variables in the $n$th Wiener chaos such that
$\lim_{N\rightarrow\infty}\mathbf{E}\left[  I_{n}(f_{N})^{2}\right]  =1$. Then the
following are equivalent:

\begin{description}
\item[(i)] As $N\rightarrow\infty$, the sequence $\{I_{n}(f_{N}):\; N\geq 1\}$ converges in distribution to a standard Gaussian random variable.

\item[(ii)] For every $\tau=1,\ldots,n-1$
\[\lim_{N\rightarrow \infty} || f_{N} \otimes_{\tau} f_{N}||^{2}_{L^{2}[[0,1]^{(2n-2\tau)}]} =0 .\]
\end{description}
\end{proposition}

\subsection{Rosenblatt process and filters: definitions, notation, and chaos
representation}

The Rosenblatt process is the (non-Gaussian) Hermite process of order 2 with
Hurst index $H\in(\frac{1}{2},1)$. It is self-similar with stationary
increments, lives in the second Wiener chaos and can be represented as a
double Wiener-It\^{o} integral of the form
\begin{equation}
Z^{(H)}(t):=Z(t)=\int_{0}^{t}\int_{0}^{t}L_{t}(y_{1},y_{2})dW_{y_{1}}%
dW_{y_{2}}. \label{rosenblatt}%
\end{equation}
Here $\{W_{t},t\in\lbrack0,1]\}$ is a standard Brownian motion and
$L_{t}(y_{1},y_{2})$ is the kernel of the Rosenblatt process
\begin{equation}
L_{t}(y_{1},y_{2})=d(H)1_{[0,t]}(y_{1})1_{[0,t]}(y_{2})\int_{y_{1}\vee y_{2}%
}^{t}\frac{\partial K^{H^{^{\prime}}}}{\partial u}(u,y_{1})\frac{\partial
K^{H^{^{\prime}}}}{\partial u}(u,y_{2})du, \label{r_kernel}%
\end{equation}
where
\[
H^{\prime}=\frac{H+1}{2}\text{ and }d(H)=\frac{1}{H+1}\left(  \frac
{H}{2(2H-1)}\right)  ^{-1/2}%
\]
and $K^{H}$ is the standard kernel of fBm, defined for $s<t$ and $H\in
(\frac{1}{2},1)\ $by%
\begin{equation}
K^{H}(t,s):=c_{H}s^{\frac{1}{2}-H}\int_{s}^{t}(u-s)^{H-\frac{3}{2}}%
u^{H-\frac{1}{2}}du \label{defK}%
\end{equation}
where $c_{H}=\left(  \frac{H(2H-1)}{\beta(2-2H,H-\frac{1}{2})}\right)
^{\frac{1}{2}}$ and $\beta(\cdot,\cdot)$ is the beta function. For $t>s$, we
have the following expression for the derivative of $K^{H}$ with respect to
its first variable:
\begin{equation}
\frac{\partial K^{H}}{\partial t}(t,s):=\partial_{1}K^{H}(t,s)=c_{H}\left(
\frac{s}{t}\right)  ^{\frac{1}{2}-H}(t-s)^{H-\frac{3}{2}}. \label{dK}%
\end{equation}
The term \textit{Rosenblatt random variable} denotes any random variable which
has the same distribution as $Z(1)$. Note that this distribution depends on
$H$.

\begin{definition}
A filter $\alpha$ of length $\ell\in\mathbb{N}$ and order $p\in\mathbb{N}\setminus {0} $ is an $(\ell+1)$-dimensional vector $\alpha
=\{\alpha_{0},\;\alpha_{1},\ldots,\alpha_{\ell}\}$ such that
\begin{align*}
\sum_{q=0}^{\ell}\alpha_{q}q^{r}  &  =0,\;\;\;\;\text{for}\;\;\;0\leq r\leq
p-1,\;r\in\mathbf{Z}\\
\sum_{q=0}^{\ell}\alpha_{q}q^{p}  &  \neq0
\end{align*}
with the convention $0^{0}=1$.
\end{definition}

If we associate such a filter $\alpha$ with the Rosenblatt process we get the
filtered process $V^{\alpha}$ according to the following scheme:%
\[
V^{\alpha}\left(  \frac{i}{N}\right)  :=\sum_{q=0}^{\ell}\alpha_{q}Z\left(
\frac{i-q}{N}\right)  ,\;\;\text{for}\;i=\ell,\ldots,N-1.
\]
Some examples are the following:

\begin{enumerate}
\item For $\alpha=\{1,-1\}$
\[
V^{\alpha}\left(  \frac{i}{N}\right)  =Z\left(  \frac{i}{N}\right)  -Z\left(
\frac{i-1}{N}\right).
\]
This is a filter of length 1 and order 1.

\item For $\alpha=\{1,-2,1\}$
\[
V^{\alpha}\left(  \frac{i}{N}\right)  =Z\left(  \frac{i}{N}\right)  -2Z\left(
\frac{i-1}{N}\right)  +Z\left(  \frac{i-2}{N}\right).
\]
This is a filter of length 2 and order 2.

\item More generally, longer filters produced by finite-differencing are such
that the coefficients of the filter $\alpha$ are the binomial coefficients
with alternating signs. Therefore, borrowing the notation $\nabla$ from time
series analysis, $\nabla Z\left(  i/N\right)  =Z\left(  i/N\right)  -Z\left(
\left(  i-1\right)  /N\right)  $, we define $\nabla^{j}=\nabla\nabla^{j-1}$
and we may write the $j$th-order finite-difference-filtered process as
follows
\[
V^{\alpha_{j}}\left(  \frac{i}{N}\right)  :=\left(  \nabla^{j}Z\right)
\left(  \frac{i}{N}\right)  .
\]

\end{enumerate}

\noindent\textbf{From now on we assume the filter order is strictly greater
than 1 ($p\geq2$).}\bigskip

For such a filter $\alpha$ the quadratic variation statistic is defined as
\[
V_{N}:=\frac{1}{N-\ell}\sum_{i=\ell}^{N-1}\left[  \frac{\left\vert V^{\alpha
}\left(  \frac{i}{N}\right)  \right\vert ^{2}}{\mathbf{E}\left\vert V^{\alpha
}\left(  \frac{i}{N}\right)  \right\vert ^{2}}-1\right]  .
\]
Using the definition of the filter, we can compute the covariance of the
filtered process $V^{\alpha}\left(  \frac{i}{N}\right)  $:%
\begin{align*}
\pi_{H}^{\alpha}(j)  &  :=\mathbf{E}\left[  V^{\alpha}\left(  \frac{i}%
{N}\right)  V^{\alpha}\left(  \frac{i+j}{N}\right)  \right] \\
&  =\sum_{q,r=0}^{\ell}\alpha_{q}\alpha_{r}\mathbf{E}\left[  Z\left(
\frac{i-q}{N}\right)  Z\left(  \frac{i+j-r}{N}\right)  \right] \\
&  =\frac{N^{-2H}}{2}\sum_{q,r=0}^{\ell}\alpha_{q}\alpha_{r}\left(
|i-q|^{2H}+|i+j-r|^{2H}-|j+q-r|^{2H}\right) \\
&  =-\frac{N^{-2H}}{2}\sum_{q,r=0}^{\ell}\alpha_{q}\alpha_{r}|j+q-r|^{2H}%
+\frac{N^{-2H}}{2}\sum_{q,r=0}^{\ell}\alpha_{q}\alpha_{r}\left(
|i-q|^{2H}+|i+j-r|^{2H}\right)  .
\end{align*}
Since the term $\sum_{q,r=0}^{\ell}\alpha_{q}\alpha_{r}\left(  |i-q|^{2H}%
+|i+j-r|^{2H}\right)  $ vanishes we get that
\begin{equation}
\pi_{H}^{\alpha}(j)=-\frac{N^{-2H}}{2}\sum_{q,r=0}^{\ell}\alpha_{q}\alpha
_{r}|j+q-r|^{2H}. \label{piH}%
\end{equation}
Therefore, we can rewrite the variation statistic as follows
\begin{align*}
V_{N}  &  =\frac{1}{N-\ell}\sum_{i=\ell}^{N-1}\left[  \frac{\left\vert
V^{\alpha}\left(  \frac{i}{N}\right)  \right\vert ^{2}}{\pi_{H}^{\alpha}%
(0)}-1\right] \\
&  =\frac{2N^{2H}}{N-\ell}\left(  -\sum_{q,r=0}^{\ell}\alpha_{r}\alpha
_{q}|q-r|^{2H}\right)  ^{-1}\sum_{i=\ell}^{N-1}\left[  \left\vert V^{\alpha
}\left(  \frac{i}{N}\right)  \right\vert ^{2}-\pi_{H}^{\alpha}(0)\right] \\
&  =\frac{2N^{2H}}{c(H)(N-\ell)}\;\sum_{i=\ell}^{N-1}\left[  \left\vert
V^{\alpha}\left(  \frac{i}{N}\right)  \right\vert ^{2}-\pi_{H}^{\alpha
}(0)\right]  ,
\end{align*}
where
\begin{equation}
c(H)=-\sum_{q,r=0}^{\ell}\alpha_{r}\alpha_{q}|q-r|^{2H}. \label{cH}%
\end{equation}
The next lemma is informative, and will be useful in the sequel.

\begin{lemma}
\label{lemch}$c\left(  H\right)  $ is positive for all $H\in(0,1]$. Also,
$c\left(  0\right)  =0$.
\end{lemma}

\begin{proof}
For $H<1$, we may rewrite $c\left(  H\right)  $ by using the representation of
the function $|q-r|^{2H}$ via fBm $B^{H}$, as its canonical metric given in
(\ref{canon}), and its covariance function $R_{H}$ given in (\ref{cov}).
Indeed we have%
\begin{align*}
c\left(  x\right)   &  =-\sum_{q,r=0}^{\ell}\alpha_{r}\alpha_{q}%
\mathbf{E}\left[  \left(  B^{H}\left(  q\right)  -B^{H}\left(  r\right)
\right)  ^{2}\right] \\
&  =-\sum_{q,r=0}^{\ell}\alpha_{r}\alpha_{q}\left(  R_{H}\left(  q,q\right)
+R_{H}\left(  r,r\right)  -2R_{H}\left(  q,r\right)  \right) \\
&  =-2\left(  \sum_{q=0}^{\ell}\alpha_{q}\right)  \left(  \sum_{r=0}^{\ell
}\alpha_{r}R_{H}\left(  r,r\right)  \right)  +2\sum_{q,r=0}^{\ell}\alpha
_{r}\alpha_{q}R_{H}\left(  q,r\right) \\
&  =0+2\sum_{q,r=0}^{\ell}\alpha_{r}\alpha_{q}R_{H}\left(  q,r\right)
=\mathbf{E}\left[  \left(  \sum_{q=0}^{\ell}\alpha_{q}B^{H}\left(  q\right)
\right)  ^{2}\right]  >0
\end{align*}
where in the second-to-last line we used the filter property which implies
$\sum_{q=0}^{\ell}\alpha_{q}=0$, and the last inequality follows from the fact
that $\sum_{q=0}^{\ell}\alpha_{q}B^{H}\left(  q\right)  $ is Gaussian and
non-constant. When $H=1$, the same argument as above holds because the
Gaussian process $X$ such that $X\left(  0\right)  =0$ and $\mathbf{E}\left[
\left(  X\left(  t\right)  -X\left(  s\right)  \right)  ^{2}\right]
=\left\vert t-s\right\vert ^{2}$ is evidently equal in law to $X\left(
t\right)  =tN$ where $N$ is a fixed standard normal r.v. The assertion that
$c(0)=0$ comes from the filter property.
\end{proof}

Observe that we can write the filtered process as an integral belonging to the
second Wiener chaos
\[
V^{\alpha}\left(  \frac{i}{N}\right)  =\sum_{q=0}^{\ell}\alpha_{q}Z\left(
\frac{i-q}{N}\right)  =I_{2}\left(  \sum_{q=0}^{\ell}\alpha_{q}L_{\frac
{i-q}{N}}\right)  :=I_{2}\left(  C_{i}\right)  ,
\]
where
\begin{equation}
C_{i}:=\sum_{q=0}^{\ell}\alpha_{q}L_{\frac{i-q}{N}}. \label{SeeEye}%
\end{equation}
Using the product formula (\ref{prod}) for multiple stochastic integrals now
results in the Wiener chaos expansion of $V_{N}$.

\begin{proposition}
\label{VNchaos}With $C_{i}$ as in (\ref{SeeEye}), the variation statistic
$V_{N}$ is given by%
\begin{align*}
V_{N}  &  =\frac{2N^{2H}}{c(H)(N-l)}\sum_{i=\ell}^{N-1}\bigg[\left\vert
I_{2}(C_{i})\right\vert ^{2}-\pi_{H}^{\alpha}(0)\bigg]\\
&  =\frac{2N^{2H}}{c(H)(N-\ell)}\left[  \sum_{i=\ell}^{N-1}I_{4}\left(
C_{i}\otimes C_{i}\right)  +4\sum_{i=\ell}^{N-1}I_{2}\left(  C_{i}\otimes
_{1}C_{i}\right)  \right] \\
&  :=T_{4}+T_{2},
\end{align*}
where $T_{4}$ is a term belonging to the 4th Wiener chaos and $T_{2}$ a term
living in the 2nd Wiener chaos.
\end{proposition}

In order to prove that a variation statistic has a normal limit we may use the
characterization of ${\mathcal{N}}(0,1)$ by Nualart and Ortiz-Latorre
(Proposition \ref{t2}). Thus, we need to start by calculating $\mathbf{E}%
\left[  |V_{N}|^{2}\right]  $ so that we can then scale appropriately, in an
attempt to apply the said proposition.

\section{Scale constants for $T_{2}$ and $T_{4}$}

In order to determine the convergence of $V_{N}$, using the orthogonality of
the integrals belonging in different chaoses, we will study each term
separately. This section begins by calculating the second moments of $T_{2}$
and $T_{4}$.

In this section we use an alternative expression for the filtered process.
More specifically, denoting $b_{q}:=\sum_{r=0}^{q}\alpha_{r}$, we rewrite
$C_{i}$ as follows, for any $i=\ell,\ldots,N-1$:%
\begin{align}
C_{i,\ell}  &  :=C_{i}=\sum_{q=0}^{\ell}\alpha_{q}L_{\frac{i-q}{N}}\nonumber\\
&  =\alpha_{0}\left(  L_{\frac{i}{N}}-L_{\frac{i-1}{N}}\right)  +(\alpha
_{0}+\alpha_{1})\left(  L_{\frac{i-1}{N}}-L_{\frac{i-2}{N}}\right)
+\ldots+(\alpha_{0}+\ldots+\alpha_{\ell-1})\left(  L_{\frac{i-(\ell-1)}{N}%
}-L_{\frac{i-\ell}{N}}\right) \nonumber\\
&  =\sum_{q=0}^{\ell}b_{q}\left(  L_{\frac{i-(q-1)}{N}}-L_{\frac{i-q}{N}%
}\right)  . \label{SeeNose}%
\end{align}
Recall that the filter properties imply $\sum_{q=0}^{\ell}\alpha_{q}=0$ and
$\alpha_{\ell}=-\sum_{q=0}^{\ell-1}\alpha_{q}$.


\subsection{Term $T_{2}$}

By Proposition \ref{VNchaos}, we can express $\mathbf{E}(T_{2}^{2})$ as:
\[
\mathbf{E}(T_{2}^{2})=\frac{64\;N^{4H}}{c(H)^{2}(N-\ell)^{2}}\mathbf{E}\left[
\left(  \sum_{i=\ell}^{N-1}I_{2}\left(  C_{i}\otimes_{1}C_{i}\right)  \right)
^{2}\right]  =\frac{2!\;64\;N^{4H}}{c(H)^{2}(N-\ell)^{2}}\sum_{i,j=\ell}%
^{N-1}\left\langle C_{i}\otimes_{1}C_{i},\;C_{j}\otimes_{1}C_{j}\right\rangle
_{L^{2}([0,1]^{2})}%
\]

\begin{proposition}
\label{propc2}We have
\[
\lim_{N\rightarrow\infty}\mathbf{E}\left[  \left\vert N^{1-H}\;T_{2}%
\right\vert ^{2}\right]  =c_{2,H}.
\]
where
\begin{equation}
c_{2,H}=\frac{64}{c(H)^{2}}\left(  \frac{2H-1}{H\;(H+1)^{2}}\right)  \left\{
\sum_{q,r=0}^{\ell}b_{q}b_{r}\left[  |1+q-r|^{2H^{^{\prime}}}%
+|1-q+r|^{2H^{^{\prime}}}-2|q-r|^{2H^{^{\prime}}}\right]  \right\}  ^{2}.
\label{c2}%
\end{equation}

\end{proposition}

This proposition is proved in the Appendix.

\subsection{Term $T_{4}$}

In this paragraph we estimate the second moment of $T_{4}$, the fourth chaos
term appearing in the decomposition of the variation $V_{N}$. Here the
function $\sum_{i=\ell}^{N-1}\left(  C_{i}\otimes C_{i}\right)  $ is no longer
symmetric and we need to symmetrize this kernel to calculate $T_{4}$'s second
moment. In other words, by Proposition \ref{VNchaos}, we have that
\begin{align*}
\mathbf{E}\left(  T_{4}^{2}\right)   &  =\frac{4N^{4H}}{c(H)^{2}(N-\ell)^{2}%
}\mathbf{E}\left[  \left(  \sum_{i=\ell}^{N-1}I_{4}(C_{i}\otimes
C_{i})\right)  ^{2}\right] \\
&  =\frac{4N^{4H}}{c(H)^{2}(N-\ell)^{2}}4!\sum_{i,j=\ell}^{N-1}\langle
C_{i}\tilde{\otimes}C_{i},C_{j}\tilde{\otimes}C_{j}\rangle_{L^{2}([0,1]^{4})}%
\end{align*}
where $C_{i}\tilde{\otimes}C_{i}:=\widetilde{C_{i}\otimes C_{i}}$. Thus, we
can use the following combinatorial formula:\newline If $f$ and $g$ are two
symmetric functions in $L^{2}([0,1]^{2})$, then
\begin{align*}
&  4!\langle f\tilde{\otimes}f,g\tilde{\otimes}g\rangle_{L^{2}([0,1]^{4})}\\
&  =(2!)^{2}\langle f\otimes f,g\otimes g\rangle_{L^{2}([0,1]^{4})}%
+(2!)^{2}\langle f\otimes_{1}g,g\otimes_{1}f\rangle_{L^{2}([0,1]^{2})}.
\end{align*}
It implies
\begin{align*}
\mathbf{E}\left(  T_{4}^{2}\right)   &  =\frac{4N^{4H}}{c(H)^{2}(N-\ell)^{2}%
}4!\sum_{i,j=\ell}^{N-1}\langle C_{i}\tilde{\otimes}C_{i},C_{j}\tilde{\otimes
}C_{j}\rangle_{L^{2}([0,1]^{4})}\\
&  =\frac{4N^{4H}}{c(H)^{2}(N-\ell)^{2}}4\;\sum_{i,j=\ell}^{N-1}\langle
C_{i}\otimes C_{i},C_{j}\otimes C_{j}\rangle_{L^{2}([0,1]^{4})}\\
&  +\frac{4N^{4H}}{c(H)^{2}(N-\ell)^{2}}4\;\sum_{i,j=\ell}^{N-1}\langle
C_{i}\otimes_{1}C_{j},C_{j}\otimes_{1}C_{i}\rangle_{L^{2}([0,1]^{2})}\\
&  :=T_{4,(1)}+T_{4,(2)}.
\end{align*}
The proof of the next proposition, in the Appendix, shows that the two terms
$T_{4,(1)}$ and $T_{4,(2)}$ have the same order of magnitude, with only the
normalizing constant being different.

\begin{proposition}
\label{propET4}Recall the constant $c\left(  H\right)  $ defined in
(\ref{cH}). Let
\begin{align*}
\tau_{1,H}  &  :=\sum_{k=\ell}^{\infty}\sum_{q_{1},q_{2},r_{1},r_{1}=0}^{\ell
}b_{q_{1}}b_{q_{2}}b_{r_{1}}b_{r_{2}}\int_{[0,1]^{4}}dudvdu^{\prime}%
dv^{\prime}\\
&  \bigg[\left\vert u-v+k-q_{1}+r_{1}\right\vert ^{2H^{\prime}-2}\left\vert
u^{\prime}-v^{\prime}+k-q_{2}+r_{2}\right\vert ^{2H^{\prime}-2}\\
&  \left\vert u-u^{\prime}+k-q_{1}+q_{2}\right\vert ^{2H^{\prime}-2}\left\vert
v-v^{\prime}+k-r_{1}+r_{2}\right\vert ^{2H^{\prime}-2}\bigg].
\end{align*}
and%
\[
\rho_{H}^{\alpha}(k):=\frac{\sum_{q,r=0}^{\ell}\alpha_{q}\alpha_{r}\left\vert
k+q-r\right\vert ^{2H}}{c(H)}%
\]
Then we have the following asymptotic variance for $\sqrt{N}T_{4}$:%
\begin{equation}
\lim_{N\rightarrow\infty}\mathbf{E}\left[  \left\vert \sqrt{N}\;T_{4}%
\right\vert ^{2}\right]  =c_{1,H}:=4!\left(  1+\sum_{k=0}^{\infty}\left\vert
\rho_{H}^{\alpha}(k)\right\vert ^{2}\right)  +\tau_{1,H}. \label{c1}%
\end{equation}

\end{proposition}

This proposition is proved in the Appendix. Observe that in the Wiener chaos
decomposition of $V_{N}$ the leading term is the term in the second Wiener
chaos (i.e. $T_{2}$) since it is of order $N^{H-1}$, while $T_{4}$ is of the
smaller order $N^{-1/2}$. We note that, in contrast to the case of filters of lenght 1 and power 1,   the barrier $H=3/4$ does not appear anymore in the estimation of the magnitude of $T_{4}$ Thus, the asymptotic behavior of $V_{N}$ is
determined by the behavior of $T_{2}$. In other words, the previous three
propositions imply the following.

\begin{theorem}
\label{th} For all $H\in(1/2,1)$ we have that
\[
\lim_{N\rightarrow\infty}\mathbf{E}\left[  \left\vert N^{1-H}\;V_{N}%
\right\vert ^{2}\right]  =c_{2,H},
\]
where $c_{2,H}$ is defined in (\ref{c2}).
\end{theorem}

>From the practical point of view, one only needs to compute the constant
$c_{2,H}$ to find the first order asymptotics of $V_{N}$. This constant is
easily computed exactly from its formula (\ref{c2}), unlike the constant
$c_{1,H}$ in Proposition (\ref{propET4}) which can only be approximated via
its unwieldy series-integral representation given therein.

\section{Normality of the term $T_{4}$}

We study in this section the limit of the renormalized term $T_{4}$ which
lives in the fourth Wiener chaos and appears in the expression of the
variation $V_{N}$. Of course, due to Theorem \ref{th} above, this term does
not affect the first order behavior of $V_{N}$ but it is interesting from the
mathematical point of view because its limit is similar to those of the
variation based on the fractional Brownian motion (\cite{TV2}). In addition,
in Section 6, we will show that the asymptotics of $T_{4}$, and indeed the
value of $c_{1,H}$, are not purely academic. They are needed in order to
calculate the asymptotic variance of the adjusted variations, those which have
a normal limit when $H\in(1/2,2/3)$.

Define the quantity
\begin{eqnarray}
G_{N}&:=&\frac{\sqrt{N}}{c_{1,H}}T_{4}=\frac{\sqrt{N}}{\sqrt{c_{1,H}}}%
\frac{2N^{2H}}{c(H)(N-\ell)}\sum_{i=\ell}^{N-1}I_{4}\left(  C_{i}\otimes
C_{i}\right) \notag\\
&=& I_{4} \left(\frac{\sqrt{N} \;2\;N^{2H}}{\sqrt{c_{1,H}}\;c(H)\;(N-\ell)} \sum_{i=\ell}^{N-1} (C_{i}\otimes C_{i}) \right) := I_{4}(g_{N}).  \label{GN}
\end{eqnarray}
>From the calculations above we proved that $\lim_{N\rightarrow\infty
}\mathbf{E}(G_{N}^{2})=1$. Using the Nualart--Peccati criterion in
Proposition \ref{t2}, we can now prove that $G_{N}$ is asymptotically standard normal.

\begin{theorem}
\label{thmT4}For all $H\in(1/2,1)$ $G_{N}$ defined in (\ref{GN}) converges in
distribution to the standard normal.
\end{theorem}

\textbf{Setup of proof of Theorem \ref{thmT4}.} To prove this theorem, by
Proposition \ref{propET4} and Proposition \ref{t2}, it is sufficient to show
that for all $\tau=1,2,3$,
\[\lim_{N\rightarrow \infty} \left \|g_{N} \tilde{\otimes_{\tau}} g_{N} \right \|_{L^{2}([0,1]^{(8-2\tau)})} = 0.\]
For  $\tau=1,2,3$, this quantity can be written as
\begin{eqnarray*}
&& \lim_{N\rightarrow \infty} \left( \frac{4 N^{4H+1}}{c_{1,H}c(H)^{2}(N-\ell)^{2}} \right)^{2}
  \left \| \sum_{i,j=\ell}^{N-1} (C_{i}\otimes C_{i}) \tilde{\otimes_{\tau}} (C_{j}\otimes C_{j})\right\|^{2}_{L^{2}([0,1]^{(8-2\tau)})}\\
&\leq& \lim_{N\rightarrow \infty} \left( \frac{4 N^{4H+1}}{c_{1,H}c(H)^{2}(N-\ell)^{2}} \right)^{2}
   \left\| \sum_{i,j=\ell}^{N-1} (C_{i}\otimes C_{i}) \otimes_{\tau} (C_{j}\otimes C_{j}) \right\|^{2}_{L^{2}([0,1]^{(8-2\tau)})}\\
&=& \lim_{N\rightarrow \infty} \left( \frac{4 N^{4H+1}}{c_{1,H}c(H)^{2}(N-\ell)^{2}} \right)^{2}
\sum_{i,j,m,n=\ell}^{N-1}      \left \langle (C_{i}\otimes C_{i}) \otimes_{\tau} (C_{j}\otimes C_{j}) ,
(C_{m}\otimes C_{m}) \otimes_{\tau} (C_{n}\otimes C_{n}) \right \rangle.
\end{eqnarray*}

The Appendix can now be consulted for proof that for each $\tau =1,2,3$ this quantity
converges to $0$, establishing the theorem.\ \ \ $\square$

\section{Anormality of the $T_{2}$ term and Asymptotic Distribution of the
2-Variation}

For the asymptotic distribution of the variation statistic we
have the following proposition.

\begin{theorem}
\label{thmV}For all $H\in(1/2,1)$, both $\frac{N^{1-H}}{\sqrt{c_{2,H}}}T_{2}$
and the normalized 2-variation $\frac{N^{1-H}}{\sqrt{c_{2,H}}}V_{N}$ converge
in $L^{2}(\Omega)$ to the Rosenblatt random variable $Z(1)$.
\end{theorem}

\textbf{Setup of proof of Theorem \ref{thmV}.} The strategy for proving this
theorem is simple. First of all Proposition \ref{propET4} implies immediately
that $N^{1-H}T_{4}$ converges to zero in $L^{2}(\Omega)$. Thus if we can show
the theorem's statement about $T_{2}$, the statement about $V_{N}$ will
following immediately from Proposition \ref{VNchaos}.

Next, to show $\frac{N^{1-H}}{\sqrt{c_{2,H}}}T_{2}$ converges to the random
variable $Z\left(  1\right)  $ in $L^{2}\left(  \Omega\right)  $, recall that
$T_{2}$\ is a second-chaos random variable of the form $I_{2}(f_{N})$, where
$f_{N}(y_{1},y_{2})$ is a symmetric function in $L^{2}([0,1]^{2})$, and that
this double Wiener-It\^{o} integral is with respect to the Brownian motion $W$
used to define $Z\left(  1\right)  $, i.e. that $Z\left(  1\right)
=I_{2}\left(  L_{1}\right)  $ where $L_{1}$ is the kernel of the Rosenblatt
process at time $1$, as defined in (\ref{r_kernel}). Therefore, by the
isometry property of Wiener-It\^{o} integrals (see (\ref{isom})), it is
necessary and sufficient to show that $\frac{N^{1-H}}{\sqrt{c_{2,H}}}f_{N}$
converges in $L^{2}([0,1]^{2})$ to $L_{1}$. This is proved in the
Appendix.\ \ \ \ $\square$

\section{Normality of the adjusted variations}

In the previous section we proved that the distribution of the variation
statistic $V_{N}$ is never normal, irrespective of the order of the filter.
However, in the decomposition of $V_{N}$, there is a normal part, $T_{4}$,
which implies that if we subtract $T_{2}$ from $V_{N}$ the remaining part will
converge to a normal law. But $T_{2}$ is not observed in practice. Following
the idea of the adjusted variations in \cite{TV}, instead of $T_{2}$ we
subtract $Z(1)$ which is observed. $Z(1)$ is the value of the Rosenblatt
process at time 1. Thus, we study the convergence of the adjusted variation:
\begin{align*}
V_{N}-\frac{\sqrt{c_{2,H}}}{N^{1-H}}Z(1)  &  =V_{N}-T_{2}+T_{2}-\frac
{\sqrt{c_{2,H}}}{N^{1-H}}Z(1)\\
&  :=T_{4}+U_{2}.%
\end{align*}

In Section 4 we showed that $\frac{\sqrt{N}}{c_{1,H}}T_{4}$ converges to a
normal law. For the quantity $U_{2}$ we prove the following proposition

\begin{proposition}
For $H\in\left(  \frac{1}{2}, \frac{2}{3} \right)  $, $\sqrt{N}U_{2}$
converges in distribution to normal with mean zero and variance given by
\begin{equation}
\label{c3}c_{3,H} := c_{2,H} \sum_{k=1}^{\infty} (N-k-1) k^{2H} F\left(
\frac{1}{k} \right)  ,
\end{equation}
where $c_{2,H}$ is defined as in (\ref{c2}) and $F$ is defined as follows
\begin{align*}
F(x)  &  = d(H)^{2}\alpha(H)^{2} \sum_{q_{1}q_{2}r_{1}r_{2}=0}^{\ell}
\int_{[0,1]^{4}} du dv du^{\prime}dv^{\prime}\left|  (u-u^{\prime}+q_{2}%
-q_{1})x+1 \right|  ^{2H^{\prime}-2}\\
&  \biggl[ \frac{128\alpha(H)^{2} d(H)^{2}}{c_{2,H} c(H)^{2}} \left|
u-v-q_{1}+r_{1} \right|  ^{2H^{\prime}-2} \left|  u^{\prime}-v^{\prime}%
-q_{2}+r_{2} \right|  ^{2H^{\prime}-2}\\
&  \left|  (v-v^{\prime}-r_{1} +r_{2})x+1 \right|  ^{2H^{\prime}-2} - \frac{16
d(H)\alpha(H)}{\sqrt{c_{2,H}} c(H)} \left|  u-v-q_{1}+r_{1} \right|
^{2H^{\prime}-2}\\
&  \left|  (v-u^{\prime}-q_{2}+r_{1})x +1\right|  ^{2H^{\prime}-2} + \left|
(u-u^{\prime}+q_{1}-q_{2})x +1 \right|  ^{2H^{\prime}-2} \biggr].
\end{align*}

\end{proposition}

\begin{proof}
The proof follows the proof of \cite[Proposition 5]{TV} and is omitted here.
\end{proof}

~

Therefore, for the adjusted variation we can prove the following

\begin{theorem}
Let $Z_{t}:t\in(0,1)$ be a Rosenblatt process with $H\in(1/2,2/3)$. Then the
adjusted variation
\[
\frac{\sqrt{N}}{c_{1,H}+c_{3,H}}\left(  V_{N}(2,\alpha)-\frac{c_{2,H}}%
{N^{1-H}}Z(1)\right).
\]
converges to a standard normal law. Here $c_{1,H},$ $c_{2,H},$ and $c_{3,H}$
are given in (\ref{c1}), (\ref{c2}), and (\ref{c3}).
\end{theorem}

\begin{proof}
The proof follows the steps of the proof of \cite[Theorem 6]{TV} and is omitted.
\end{proof}

\section{Estimators for the self-similarity index}

We construct estimators for the self-similarity index of a Rosenblatt process
$Z$ based on the discrete observations at times $0,\frac{1}{N},\frac{2}%
{N},\ldots,1$. Their strong consistency and asymptotic distribution will be
consequences of the theorems above.

\subsection{Setup of the estimation problem}

Consider the quadratic variation statistic for a filter $\alpha$ of order $p$
based on the observations of our Rosenblatt process $Z$:
\begin{equation}
S_{N}:=\frac{1}{N}\sum_{i=\ell}^{N}\left(  \sum_{q=0}^{\ell}\alpha_{q}Z\left(
\frac{i-q}{N}\right)  \right)  ^{2}. \label{SN}%
\end{equation}
We have already established that $\mathbf{E}\left[  S_{N}\right]
=-\frac{N^{-2H}}{2}\sum_{q,r=0}^{\ell}\alpha_{q}\alpha_{r}|q-r|^{2H}$ (see
expression (\ref{piH}) ). By considering that $\mathbf{E}\left[  S_{N}\right]
$ can be estimated by the empirical value $S_{N}$, we can construct an
estimator $\hat{H}_{N}$ for $H$ by solving the following equation:%
\[
S_{N}=-\frac{N^{-2\hat{H}_{N}}}{2}\sum_{q,r=0}^{\ell}\alpha_{q}\alpha
_{r}|q-r|^{2\hat{H}_{N}}.%
\]
In this case, unlike the case of a filter of length $1$ which was studied in
\cite{TV}, we cannot compute an analytical expression for the estimator.
Nonetheless, the estimator $\hat{H}_{N}$ can be easily computed numerically by
solving the following non-linear equation for fixed $N$, with unknown
$x\in\lbrack1/2,1]$:%
\begin{equation}
-\frac{N^{-2x}}{2}\sum_{q,r=0}^{\ell}\alpha_{q}\alpha_{r}|q-r|^{2x}%
-S_{N}(2,\alpha)=0. \label{hhateq}%
\end{equation}
This equation is not entirely trivial, in the sense that one must determine
whether it has a solution in $[1/2,1]$, and whether this solution is unique.
As it turns out, the answer to both questions is affirmative for large $N$, as
seen in the next Proposition, proved further below.

\begin{proposition}
\label{propex}Almost surely, for large $N$, equation (\ref{hhateq}) has
exactly one solution in $[1/2,1]$.
\end{proposition}

\begin{definition}
We define the estimator $\hat{H}_{N}$ of $H$ to be the unique solution of
(\ref{hhateq}).
\end{definition}

Note that Equation (\ref{hhateq}) can be rewritten as $S_{N}=c(x)N^{-2x}/2$
where the function $c$ was defined in (\ref{cH}). The proposition is
established via the following lemma.

\begin{lemma}
\label{lemNSN}For any $H\in(1/2,1)$, almost surely, $\lim_{N\rightarrow\infty
}N^{2H}S_{N}=c\left(  H\right)  /2$.
\end{lemma}

\begin{proof}
Firstly, we show that $V_{N}$ converges to zero almost surely as
$N\rightarrow\infty$. We already know that this is true in $L^{2}\left(
\Omega\right)  $. Consider the following
\[
P\left(  |V_{N}|>N^{-\beta}\right)  \leq N^{q\beta}\mathbf{E}\left(
|V_{N}|^{q}\right)  \leq c_{q,4}\left[  \mathbf{E}\left(  V_{N}^{2}\right)
\right]  ^{q/2}\leq c\;N^{q\beta}N^{(H-1)q}.
\]
If we choose $\beta<1-H$ and $q$ large enough so that $(1-H-\beta)q>1$. This
implies that
\[
\sum_{N=0}^{\infty}P\left(  |V_{N}|>N^{-\beta}\right)  \leq c\sum
_{N=0}^{\infty}N^{(\beta+H-1)q}<+\infty
\]
Therefore, the Borel-Cantelli lemma implies $|V_{N}|\rightarrow0$ a.s., with
speed of convergence equal to $N^{-\beta}$, for all $\beta<1-H$. Since
$V_{N}=\frac{S_{N}}{\mathbf{E}(S_{N})}-1$ we have%
\begin{equation}
1+V_{N}=-\frac{2N^{2H}}{\sum_{q,r=0}^{\ell}\alpha_{q}\alpha_{r}|q-r|^{2H}%
}S_{N}=2N^{2H}S_{N}/c\left(  H\right)  . \label{lemNSNline}%
\end{equation}
The almost-sure convergence of $V_{N}$ to $0$ yields the statement of the lemma.
\end{proof}

\bigskip

\begin{proof}
[Proof of Proposition \ref{propex}]For $x\in\lbrack\frac{1}{2},1]$ and for any
fixed $N$, define the function
\[
F_{N}\left(  x\right)  =\frac{c(x)}{2}N^{-2x}-S_{N}=-\frac{N^{-2x}}{2}%
\sum_{q,r=0}^{\ell}\alpha_{q}\alpha_{r}|q-r|^{2x}-S_{N}.%
\]
Equation (\ref{hhateq}) is $F_{N}\left(  x\right)  =0$. Observe that
$F_{N}(x)$ is strictly decreasing. Indeed, we have that
\[
F_{N}^{\prime}\left(  x\right)  =\log\left(  N^{-2x}\right)  \;\sum
_{q,r=0}^{\ell}\alpha_{q}\alpha_{r}|q-r|^{2x}-N^{-2x}\sum_{q,r=0}^{\ell}%
\alpha_{q}\alpha_{r}\log|q-r|\;|q-r|^{2x}.%
\]
Then, $F_{N}^{\prime}\left(  x\right)  <0$ is equivalent to
\[
N>\exp\left\{  \frac{\sum_{q,r=0}^{\ell}\alpha_{q}\alpha_{r}\log
|q-r|\;|q-r|^{2x}}{\sum_{q,r=0}^{\ell}\alpha_{q}\alpha_{r}|q-r|^{2x}}\right\}
,
\]
since we know, using Lemma \ref{lemch}, that $c\left(  x\right)  =\sum
_{q,r=0}^{\ell}\alpha_{q}\alpha_{r}|q-r|^{2x}$, which is evidently continuous
on $[\frac{1}{2},1]$, is strictly negative on that interval. Thus, if we
choose $N$ to be large enough, i.e.
\[
N>\max_{x\in\lbrack\frac{1}{2}{1}]}\exp\left\{  \frac{\sum_{q,r=0}^{\ell
}\alpha_{q}\alpha_{r}\log|q-r|\;|q-r|^{2x}}{\sum_{q,r=0}^{\ell}\alpha
_{q}\alpha_{r}|q-r|^{2x}}\right\}
\]
the function $F_{N}\left(  x\right)  $ is invertible on $[\frac{1}{2},1]$, and
equation (\ref{hhateq}) has no more than one solution there.

To guarantee existence of a solution, we use Lemma \ref{lemNSN}. This lemma
implies the existence of a sequence $\varepsilon_{N}$ such that $2N^{2H}%
S_{N}=c(H)+\varepsilon_{N}$ and $\lim_{N\rightarrow\infty}\varepsilon_{N}=0$
almost surely. Since in addition $c$ is continuous, then almost surely, we can
choose $N$ large enough, so that $2N^{2H}S_{N}$ is in the image of $[\frac
{1}{2},1]$ by the function $c$. Thus the equation $c\left(  x\right)
=2N^{2H}S_{N}$ has at least one solution in $[\frac{1}{2},1]$. Since this
equation is equivalent to (\ref{hhateq}), the proof of the proposition is complete.
\end{proof}

\subsection{Properties of the estimator}

Now, it remains to prove that any such $\hat{H}_{N}$ is consistent and to
determine its asymptotic distribution.

\begin{theorem}
\label{thmconsistent}For $H\in(1/2,1)$ assume that the observed process used
in the previous definition is a Rosenblatt process with Hurst parameter $H$.
Then strong consistency holds for $\hat{H}_{N}$, i.e.
\[
\lim_{N\rightarrow\infty}\hat{H}_{N}=H,\;\;a.s.
\]
In fact, we have more precisely that $\lim_{N\rightarrow\infty}\left(
H-\hat{H}_{N}\right)  \log N=0$ a.s.
\end{theorem}

\begin{proof}
>From line (\ref{lemNSNline}) in the proof of Lemma \ref{lemNSN}, and using the
fact that $\hat{H}_{N}$ solves equation (\ref{hhateq}), i.e. $c\left(  \hat
{H}_{N}\right)  N^{-\hat{H}_{N}}=2S_{N}$, we can write%
\[
1+V_{N}=-\frac{2N^{2H}}{\sum_{q,r=0}^{\ell}\alpha_{q}\alpha_{r}|q-r|^{2H}%
}S_{N}=\frac{c(\hat{H}_{N})}{c\left(  H\right)  }N^{2(H-\hat{H}_{N})}.
\]

Now note that $c(\hat{H}_{N})/c\left(  H\right)  $ is the ratio of two values
of the continuous function $c$ at two points in $[1/2,1]$. However, Lemma
\ref{lemch} proves that on this interval, the function $c$ is strictly
positive; since it is continuous, it is bounded above and away from $0$. Let
$a=\min_{x\in\lbrack1/2,1]}c\left(  x\right)  >0$ and $A=\max_{x\in
\lbrack1/2,1]}c\left(  x\right)  <\infty$. These constants $a$ and $A$ are of
course non random. Therefore $c(\hat{H}_{N})/c\left(  H\right)  $ is always in
the interval $[a/A,A/a]$. Thus, almost surely,
\[
\left\vert \log\left(  c(\hat{H}_{N})/c\left(  H\right)  \right)  \right\vert
\leq\log\frac{A}{a}.
\]
We may now write%
\begin{equation}
\log\left(  1+V_{N}\right)  =2\left(  H-\hat{H}_{N}\right)  \log N+\log\left(
\frac{c(\hat{H}_{N})}{c\left(  H\right)  }\right)  . \label{thmconsistentline}%
\end{equation}
Since in addition $\lim_{N\rightarrow\infty}\log\left(  1+V_{N}\right)  =0$
a.s., we get that almost surely,
\[
\left\vert H-\hat{H}_{N}\right\vert =O\left(  \frac{1}{\log N}\right)  .
\]
This implies the first statement of the proposition.

The second statement, which is more precise, is now obtained as follows. Since
$\hat{H}_{N}\rightarrow H$ almost surely, and $c$ is continuous, $\log\left(
c(\hat{H}_{N})/c\left(  H\right)  \right)  $ converges to $0$. The second
statement now follows immediately.
\end{proof}

\bigskip

The asymptotic distribution of the estimator $\hat{H}_{N}$ is stated in the
next result. Its proof uses Theorem \ref{thmV} and Theorem \ref{th}, plus the
expression (\ref{thmconsistentline}). While novel and interesting, this proof
is more technical than the proofs of the proposition and theorem above, and is
therefore relegated to the Appendix.

\begin{theorem}
\label{thmstatconv}For any $H\in(\frac{1}{2},1)$, the convergence%
\[
\lim_{N\rightarrow\infty}2c_{2,H}^{-1/2}N^{1-H}\left(  \hat{H}_{N}-H\right)
\log N=Z(1)
\]
holds in $L^{2}\left(  \Omega\right)  $, where $Z(1)$ is a Rosenblatt random variable.
\end{theorem}

As can be seen from Theorem \ref{thmV} and Theorem \ref{thmstatconv}, the
renormalization of the statistic $V_{N}$, as well as the renormalization of
the difference $\hat{H}_{N}-H$, depend on $H$: it is of order of $N^{1-H}$.
The quantities $N^{1-H}V_{N}$ and $N^{1-H}\hat{H}_{N}$ cannot be computed
numerically from the empirical data, thereby compromising the use of the
asymptotic distributions for statistical purposes such as model validation.
Therefore one would like to have other quantities with known asymptotic
distribution which can be calculated using only the data. The next theorem
addresses this issue by showing that one can replace $H$ by $\hat{H}_{N}$ in
the term $N^{1-H}$, and still obtain a convergence as in Theorem
\ref{thmstatconv}, this time in $L^{1}\left(  \Omega\right)  $. Its proof is
in the Appendix.

\begin{theorem}
\label{ThmStat}For any $H\in(\frac{1}{2},1)$, with the Rosenblatt random
variable $Z\left(  1\right)  $,%
\[
\lim_{N\rightarrow\infty}\mathbf{E}\left[  \left\vert 2\;c_{2,H}%
^{-1/2}N^{1-\hat{H}_{N}}\log N\left(  \hat{H}_{N}-H\right)  -Z\left(
1\right)  \right\vert \right]  =0.
\]

\end{theorem}

\section{Numerical Computation of the Asymptotic Variance}

In practice certain issues may occur when we compute the asymptotic variance.
The most crucial question is what order of filter we should choose. Indeed,
from (\ref{c2}) with $\hat{H}_{N}$ instead of $H$, it follows that the
constants of the variance not only depend on the filter length/order ($\ell$,
$p$), but also on the number of observations ($N$). We measure the
\textquotedblleft accuracy\textquotedblright\ of the estimator $\hat{H}_{N}$
by its standard error which is the following quantity:
\[
\frac{\sqrt{c_{2,\hat{H}_{N}}}}{2N^{1-\hat{H}_{N}}\log N}.
\]
There are several types of filters that we can use. In this paper, we choose
to work with finite-difference and wavelet-type filters.

\begin{itemize}
\item The finite-difference filters are produced by finite-differencing the
process. In this case the filter length is the same as the order of the
filter. The coefficients of the order-$\ell$ finite difference filter are
given by
\[
\alpha_{k}=(-1)^{k+1}\binom{\ell}{k},\;\;k=0,\ldots,\ell.
\]

\item The wavelet filters we are using are the Daubechies filters with
$k$-vanishing moments. (By vanishing moments we mean that all moments of the
wavelet filter are zero up to a power). The Daubechies wavelets form a family
of orthonormal wavelets with compact support and the maximum number of
vanishing moments. In this scenario, the number of vanishing moments
determines the order of the filter and the filter length is twice the order.
For more details, the reader can refer to \cite{LMR}.
\end{itemize}

We computed the standard error for $N=10,000$ observations, filters of order
varying from 2 to 20 and Hurst parameters varying from 0.55 to 0.95. This
means that the corresponding lengths of the finite-difference filters were 2
to 20 and for the wavelets 4 to 40. The code we use to simulate the Rosenblatt
process is based on a Donsker-type limit theorem and was provided to us by
J.M. Bardet \cite{Ba}. The results are illustrated in the figures 1, 2, and 3,
on the next page; these are graphs of the asymptotic standard error
$\sqrt{c_{2,H}}/(2N^{1-\hat{H}_{N}}\log N)$ for various fixed values of $H$ as
the order of the filters increase.

We observe that the standard error decreases with the order of the filter.
Furthermore, we observe that the wavelet filters are more effective than the
finite-difference ones, since they have a higher impact on the decrease of the
standard error for the same order, as the filter increases. Specifically, the
graph in Fig. 1, with the finite difference filters, shows that for fixed $H$,
there is no advantage to using a filter beyond a certain order $p$, since the
standard error tends to a constant as $p\rightarrow\infty$. This does not
occur for the wavelet filters, where the standard error continues to decrease
as $p\rightarrow\infty$ in all cases as seen in the graph in Fig. 2. On the
other hand, the finite-difference filters have lower errors than the wavelet
filters for low filter lengths; only after a certain order $p^{\ast}$ do the
latter become more effective; this comparison is seen in the graph in Fig. 3,
where $p^{\ast}$ is roughly 9.

In addition, since the order of convergence depends on the true value of the
Hurst parameter $H$, we investigated the behavior of the error with respect to
$H$. It seems that the higher $H$ is, the more we lose in terms of accuracy;
this is visible in all three graphs.

In general, the choice of a longer filter might lead to a smaller error, but
at the same time it increases the computational time needed in order to
compute $\hat{H}$ and its standard error. In a future work, we will study
extensively this trade-off and other consequences of using longer filters.

\begin{figure}[ptb]
\centering%
\begin{tabular}
[c]{cc}%
\includegraphics [width=0.4\linewidth, height=0.5\linewidth]{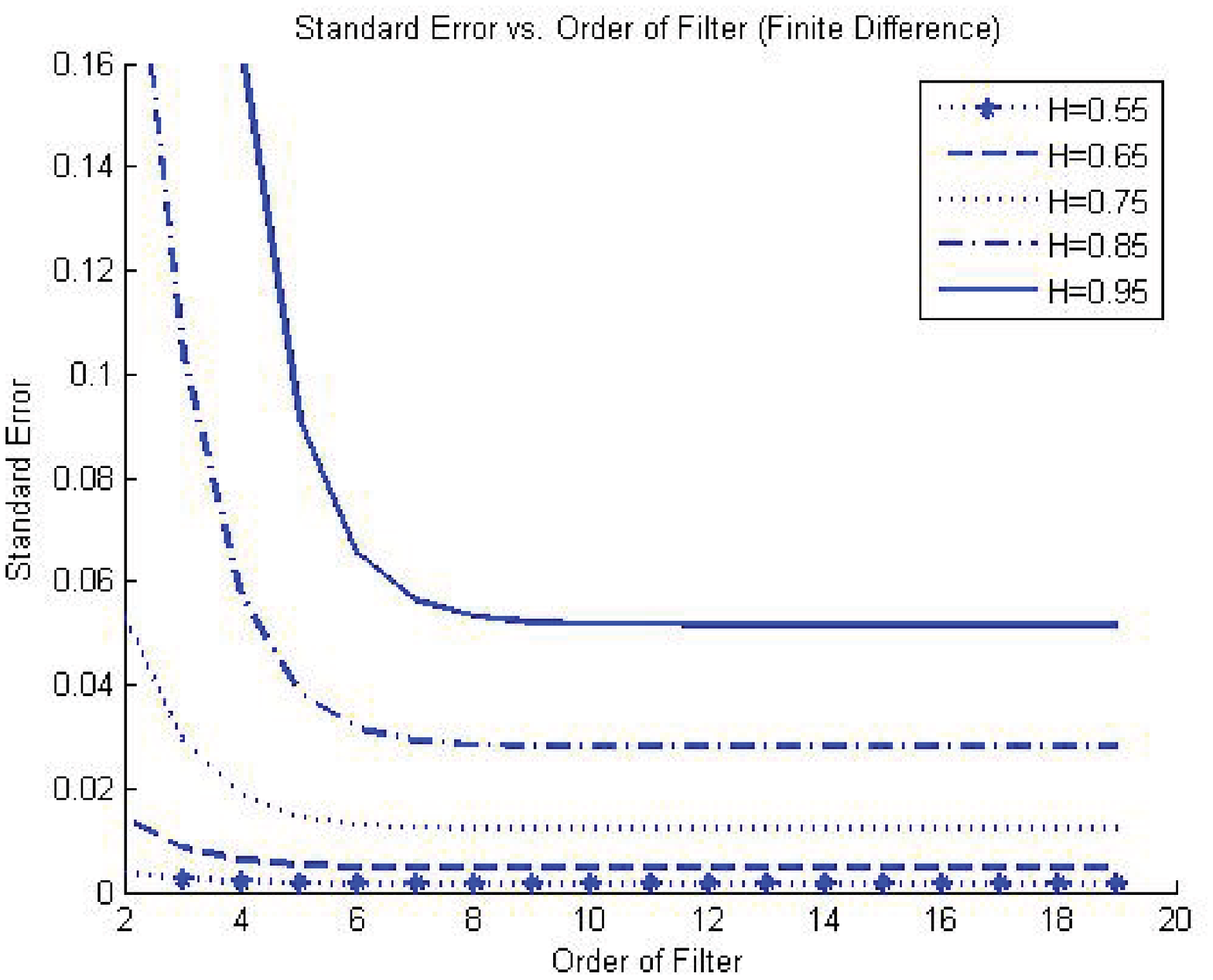} &
\includegraphics [width=0.4\linewidth, height=0.5\linewidth]{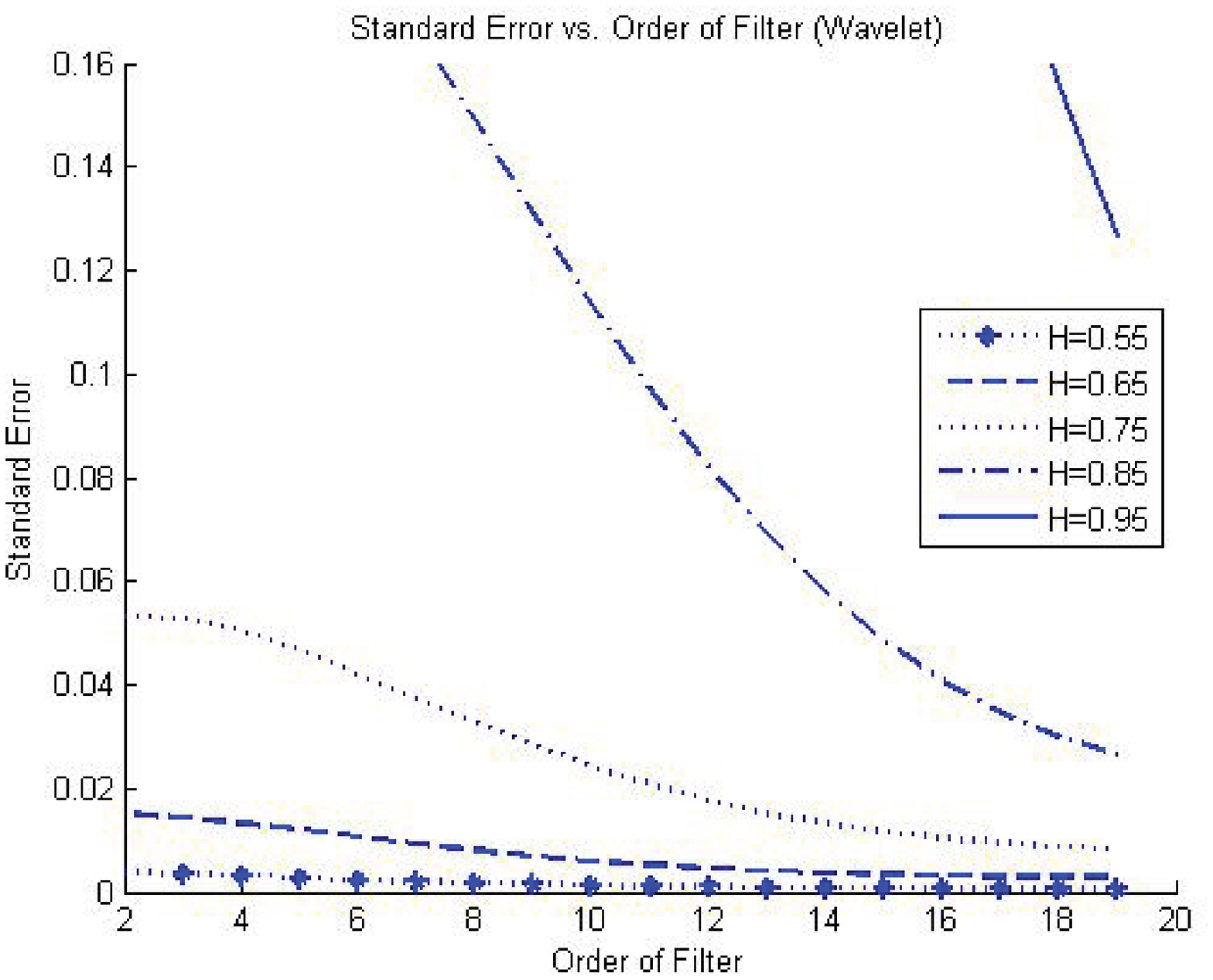}\\
Fig. 1: Finite Difference Filters. & Fig. 2.: Wavelet Filters.\\
\vspace{1cm} &
\end{tabular}
\centering
\includegraphics[width=0.5\linewidth, height=0.5\linewidth]{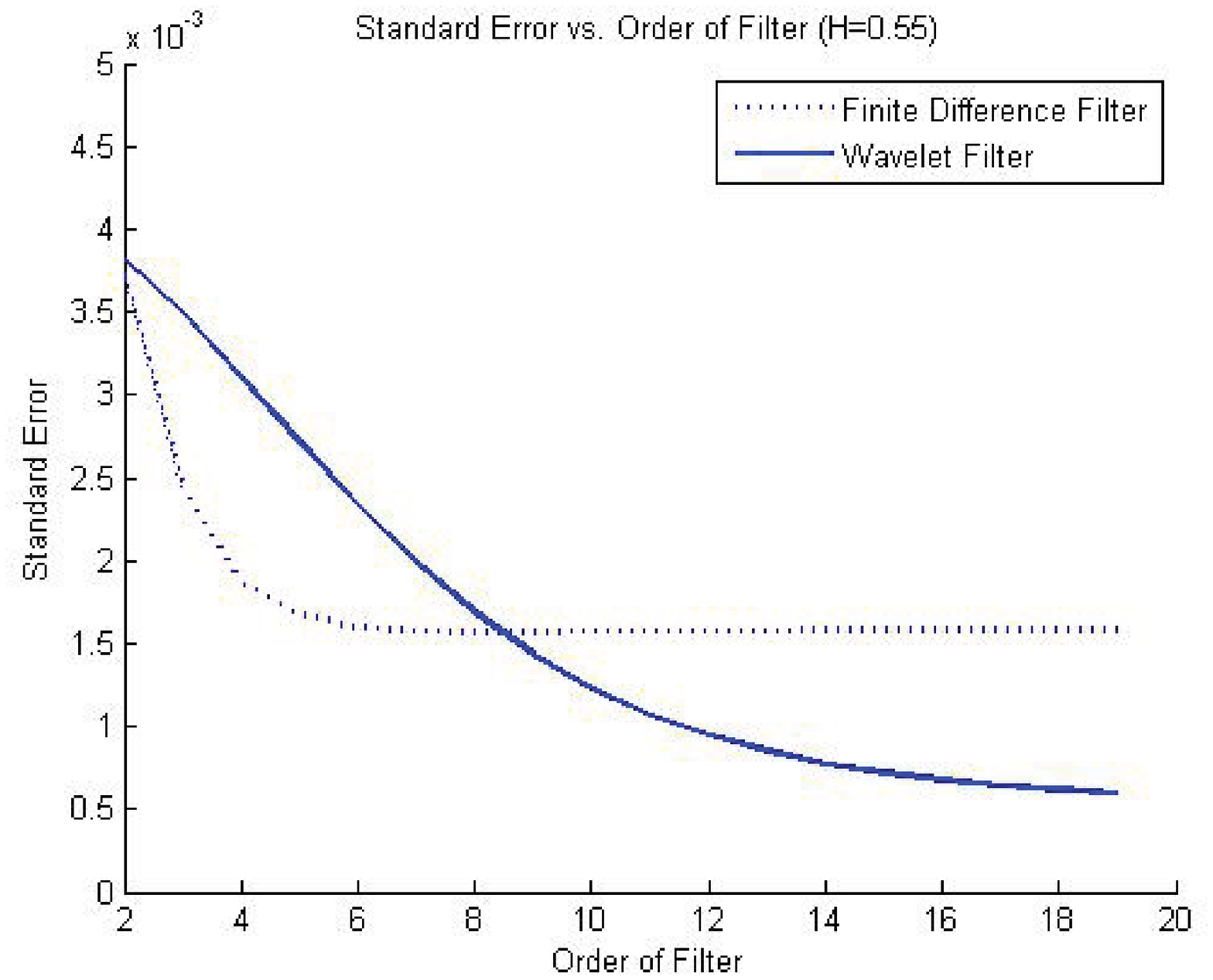}
\newline Fig. 3.: Comparison between the two types of
filter.\end{figure}\label{graph}

\bigskip

\section{Appendix: proofs.}

\subsection{Proof of Proposition \ref{propc2}}

We start by computing the contraction term $C_{i}\otimes_{1}C_{i}$:%
\begin{align*}
&  (C_{i}\otimes_{1}C_{i})(y_{1},y_{2})=\int_{0}^{1}C_{i}(x,y_{1}%
)C_{i}(x,y_{2})dx\\
&  =\sum_{q,r=0}^{\ell}b_{q}b_{r}\int_{0}^{1}\left(  L_{\frac{i-(q-1)}{N}%
}(x,y_{1})-L_{\frac{i-q}{N}}(x,y_{1})\right)  \left(  L_{\frac{i-(r-1)}{N}%
}(x,y_{2})-L_{\frac{i-r}{N}}(x,y_{2})\right)  dx\\
&  =d(H)^{2}\sum_{q,r=0}^{\ell}b_{q}b_{r}1_{[0,\frac{i-q+1}{N}]}%
(y_{1})1_{[0,\frac{i-r+1}{N}]}(y_{2})\int_{0}^{\frac{i-q+1}{N}\wedge
\frac{i-r+1}{N}}dx\\
&  \times\left(  \int_{\frac{i-q}{N}}^{\frac{i-q+1}{N}}\frac{\partial
K^{H^{^{\prime}}}}{\partial u}(u,x)\frac{\partial K^{H^{^{\prime}}}}{\partial
u}(u,y_{1})du\right)  \left(  \int_{\frac{i-r}{N}}^{\frac{i-r+1}{N}}%
\frac{\partial K^{H^{^{\prime}}}}{\partial v}(v,x)\frac{\partial
K^{H^{^{\prime}}}}{\partial v}(v,y_{2})dv\right) \\
&  =d(H)^{2}\sum_{q,r=0}^{\ell}b_{q}b_{r}1_{[0,\frac{i-q+1}{N}]}%
(y_{1})1_{[0,\frac{i-r+1}{N}]}(y_{2})\\
&  \times\int_{I_{i_{q}}}\int_{I_{i_{r}}}du\;dv\frac{\partial K^{H^{^{\prime}%
}}}{\partial u}(u,y_{1})\frac{\partial K^{H^{^{\prime}}}}{\partial u}%
(v,y_{2})dudv\left(  \int_{0}^{u\wedge v}dx\frac{\partial K^{H^{^{\prime}}}%
}{\partial u}(u,x)\frac{\partial K^{H^{^{\prime}}}}{\partial v}(v,x)\right) \\
&  =\alpha(H)d(H)^{2}\sum_{q,r=0}^{\ell}b_{q}b_{r}1_{[0,\frac{i-q+1}{N}%
]}(y_{1})1_{[0,\frac{i-r+1}{N}]}(y_{2})\\
& \int_{I_{i_{q}}}\int_{I_{i_{r}}}du\;dv|u-v|^{2H^{^{\prime}}-2}\frac{\partial K^{H^{^{\prime}}}}{\partial
u}(u,y_{1})\frac{\partial K^{H^{^{\prime}}}}{\partial v}(v,y_{2})dudv,
\end{align*}
where $I_{i_{q}}=\left(  \frac{i-q}{N},\frac{i-q+1}{N}\right]  $.

Now, the inner product computes as
\begin{align*}
&  \left\langle C_{i}\otimes_{1}C_{i},C_{j}\otimes_{1}C_{j}\right\rangle
_{L^{2}[0,1]^{2}}\\
&  =\alpha(H)^{2}d(H)^{4}\sum_{q_{1},r_{1},q_{2},r_{2}=0}^{\ell}b_{q_{1}%
}b_{r_{1}}b_{q_{2}}b_{r_{2}}\int_{0}^{1}\int_{0}^{1}dy_{1}dy_{2}\\
& \int_{I_{i_{q_{1}}}}\int_{I_{i_{r_{1}}}}\int_{I_{j_{q_{2}}}}\int_{I_{j_{r_{2}%
}}}dudvdu^{^{\prime}}dv^{^{\prime}} |u-v|^{2H^{^{\prime}}-2} |u^{^{\prime}}-v^{^{\prime}}|^{2H^{^{\prime}}-2}\\
& \frac{\partial K^{H^{^{\prime}}}}{\partial u}(u,y_{1})\frac{\partial K^{H^{^{\prime}}}%
}{\partial v}(v,y_{2})\frac{\partial K^{H^{^{\prime}}}}{\partial u^{^{\prime}%
}}(u^{^{\prime}},y_{1})\frac{\partial K^{H^{^{\prime}}}}{\partial v^{^{\prime
}}}(v^{^{\prime}},y_{2})dudvdu^{^{\prime}}dv^{^{\prime}}\\
&  =\alpha(H)^{2}d(H)^{4}\sum_{q_{1},r_{1},q_{2},r_{2}=0}^{\ell}b_{q_{1}%
}b_{r_{1}}b_{q_{2}}b_{r_{2}}\\
& \int_{I_{i_{q_{1}}}}\int_{I_{i_{r_{1}}}}%
\int_{I_{j_{q_{2}}}}\int_{I_{j_{r_{2}}}}dudvdu^{^{\prime}}dv^{^{\prime}%
}|u-v|^{2H^{^{\prime}}-2}|u^{^{\prime}}-v^{^{\prime}}|^{2H^{^{\prime}}-2}\\
&  \left(  \int_{0}^{u\wedge u^{\prime}}\frac{\partial K^{H^{^{\prime}}}%
}{\partial u}(u,y_{1})\frac{\partial K^{H^{^{\prime}}}}{\partial u^{^{\prime}%
}}(u^{^{\prime}},y_{1})dy_{1}\right)  \left(  \int_{0}^{v\wedge v^{\prime}%
}\frac{\partial K^{H^{^{\prime}}}}{\partial u}(u,y_{1})\frac{\partial
K^{H^{^{\prime}}}}{\partial v^{^{\prime}}}(v^{^{\prime}},y_{2})dy_{2}\right)
\\
&  =\alpha(H)^{4}d(H)^{4}\sum_{q_{1},r_{1},q_{2},r_{2}=0}^{\ell}b_{q_{1}%
}b_{r_{1}}b_{q_{2}}b_{r_{2}}\int_{I_{i_{q_{1}}}}\int_{I_{i_{r_{1}}}}%
\int_{I_{j_{q_{2}}}}\int_{I_{j_{r_{2}}}}dudvdu^{^{\prime}}dv^{^{\prime}}\\
&  \times|u-v|^{2H^{^{\prime}}-2}|u^{^{\prime}}-v^{^{\prime}}|^{2H^{^{\prime}%
}-2}|u-u^{^{\prime}}|^{2H^{^{\prime}}-2}|v-v^{^{\prime}}|^{2H^{^{\prime}}-2}.
\end{align*}
We make the following change of variables
\[
\bar{u}=\left(  u-\frac{i-q_{1}}{N}\right)  N
\]
and the second moment of $T_{2}$ becomes
\begin{align*}
&  \mathbf{E}\left(  T_{2}^{2}\right) \\
&  =\frac{128\;\alpha(H)^{4}d(H)^{4}}{c(H)^{2}}\frac{N^{4H}}{(N-\ell)^{2}}%
\sum_{i,j=\ell}^{N-1}\;\;\sum_{q_{1},r_{1},q_{2},r_{2}=0}^{\ell}b_{q_{1}%
}b_{r_{1}}b_{q_{2}}b_{r_{2}}\int_{I_{i_{q_{1}}}}\int_{I_{i_{r_{1}}}}%
\int_{I_{j_{q_{2}}}}\int_{I_{j_{r_{2}}}}dudvdu^{^{\prime}}dv^{^{\prime}}\\
&  \times|u-v|^{2H^{^{\prime}}-2}|u^{^{\prime}}-v^{^{\prime}}|^{2H^{^{\prime}%
}-2}|u-u^{^{\prime}}|^{2H^{^{\prime}}-2}|v-v^{^{\prime}}|^{2H^{^{\prime}}-2}\\
&  =\frac{128\;\alpha(H)^{4}d(H)^{4}}{c(H)^{2}}\frac{N^{4H}}{(N-\ell)^{2}%
}\frac{1}{N^{4}N^{8H^{^{\prime}}-8}}\sum_{i,j=\ell}^{N-1}\;\;\sum_{q_{1}%
,r_{1},q_{2},r_{2}=0}^{\ell}b_{q_{1}}b_{r_{1}}b_{q_{2}}b_{r_{2}}%
\int_{[0,1]^{4}}dudvdu^{^{\prime}}dv^{^{\prime}}\\
&  \times|u-v-q_{1}+r_{1}|^{2H^{^{\prime}}-2}|u^{^{\prime}}-v^{^{\prime}%
}-q_{2}+r_{2}|^{2H^{^{\prime}}-2}\\
&  \times|u-u^{^{\prime}}+i-j-q_{1}+q_{2}|^{2H^{^{\prime}}-2}|v-v^{^{\prime}%
}+i-j-r_{1}+r_{2}|^{2H^{^{\prime}}-2}\\
&  =\frac{128\;\alpha(H)^{4}d(H)^{4}}{c(H)^{2}}\frac{1}{(N-\ell)^{2}}%
\;\;\sum_{i,j=\ell}^{N-1}\sum_{q_{1},r_{1},q_{2},r_{2}=0}^{\ell}b_{q_{1}%
}b_{r_{1}}b_{q_{2}}b_{r_{2}}\int_{[0,1]^{4}}dudvdu^{^{\prime}}dv^{^{\prime}}\\
&  \times|u-v-q_{1}+r_{1}|^{2H^{^{\prime}}-2}|u^{^{\prime}}-v^{^{\prime}%
}-q_{2}+r_{2}|^{2H^{^{\prime}}-2}\\
&  \times\left(  |u-u^{^{\prime}}+i-j-q_{1}+q_{2}|^{2H^{^{\prime}}%
-2}|v-v^{^{\prime}}+i-j-r_{1}+r_{2}|^{2H^{^{\prime}}-2}\right)  .
\end{align*}
Let $cst.=\frac{128\;\alpha(H)^{4}d(H)^{4}}{c(H)^{2}}$. We study first the
diagonal terms of the above double sum
\begin{align*}
&  \mathbf{E}\left(  T_{2-diag}^{2}\right) \\
&  =cst.\frac{N-\ell-1}{(N-\ell)^{2}}\sum_{q_{1},r_{1},q_{2},r_{2}=0}^{\ell
}b_{q_{1}}b_{r_{1}}b_{q_{2}}b_{r_{2}}\int_{[0,1]^{4}}dudvdu^{^{\prime}%
}dv^{^{\prime}}\\
&  \times|u-v-q_{1}+r_{1}|^{2H^{^{\prime}}-2}|u^{^{\prime}}-v^{^{\prime}%
}-q_{2}+r_{2}|^{2H^{^{\prime}}-2}|u-u^{^{\prime}}-q_{1}+q_{2}|^{2H^{^{\prime}%
}-2}|v-v^{^{\prime}}-r_{1}+r_{2}|^{2H^{^{\prime}}-2}.
\end{align*}
We conclude that
\[
\mathbf{E}\left(  T_{2-diag}^{2}\right)  ={\mathcal{O}}\left(  N^{-1}\right)
.
\]
Let's consider now the non-diagonal terms
\begin{align}
&  \mathbf{E}\left(  T_{2-off}^{2}\right) = 2cst.\sum_{q_{1},r_{1},q_{2},r_{2}=0}^{\ell}b_{q_{1}}b_{r_{1}}b_{q_{2}%
}b_{r_{2}}\nonumber \\
& \times \int_{[0,1]^{4}}dudvdu^{^{\prime}}dv^{^{\prime}}\times
|u-v-q_{1}+r_{1}|^{2H^{^{\prime}}-2}|u^{^{\prime}}-v^{^{\prime}}-q_{2}%
+r_{2}|^{2H^{^{\prime}}-2}\nonumber\\
&  \times\frac{1}{(N-\ell)^{2}}\left(  \sum_{i,j=\ell,\;i\neq j}%
^{N-1}|u-u^{^{\prime}}+i-j-q_{1}+q_{2}|^{2H^{^{\prime}}-2}|v-v^{^{\prime}%
}+i-j-r_{1}+r_{2}|^{2H^{^{\prime}}-2}.\right)  \label{r_sum}%
\end{align}
Observe that the term (\ref{r_sum}) can be calculated as follows:
\begin{align*}
&  \frac{1}{(N-\ell)^{2}}\sum_{i,j=\ell\;i\neq j}^{N-1}|u-u^{^{\prime}%
}+i-j-q_{1}+r_{1}|^{2H^{^{\prime}}-2}|v-v^{^{\prime}}+i-j-r_{1}+r_{2}%
|^{2H^{^{\prime}}-2}\\
&  =\frac{1}{(N-\ell)^{2}}\sum_{i=\ell}^{N-1}\sum_{k=1}^{N-i}|u-u^{^{\prime}%
}+k-q_{1}+q_{2}|^{2H^{^{\prime}}-2}|v-v^{^{\prime}}+k-r_{1}+r_{2}%
|^{2H^{^{\prime}}-2}\\
&  =\frac{1}{(N-\ell)^{2}}\sum_{k=\ell}^{N-1}(N-k-1)|u-u^{^{\prime}}%
+k-q_{1}+q_{2}|^{2H^{^{\prime}}-2}|v-v^{^{\prime}}+k-r_{1}+r_{2}%
|^{2H^{^{\prime}}-2}\\
&  =N^{4H^{^{\prime}}-4}\frac{N}{(N-\ell)^{2}} \\
& \times \sum_{k=\ell}^{N-1}\left(
1-\frac{k+1}{N}\right)  \left\vert \frac{u-u^{^{\prime}}}{N}+\frac{k}{N}%
-\frac{q_{1}-q_{2}}{N}\right\vert ^{2H^{^{\prime}}-2}\left\vert \frac
{v-v^{^{\prime}}}{N}+\frac{k}{N}-\frac{r_{1}-r_{2}}{N}\right\vert
^{2H^{^{\prime}}-2}.
\end{align*}
We may now use a Riemann sum approximation and the fact that $4H^{^{\prime}%
}-4=2H-2>-1$. Since $\ell$ is fixed and $q_{1}$ and $q_{2}$ are less than
$\ell$, we get that the term in (\ref{r_sum}) is asymptotically equivalent to%
\[
\sum_{k=\ell}^{N-1}\left(  1-\frac{k}{N}\right)  \left\vert \frac{k}%
{N}\right\vert ^{2H^{^{\prime}}-2}\left\vert \frac{k}{N}\right\vert
^{2H^{^{\prime}}-2}=\int_{0}^{1}\left(  1-x\right)  x^{2H-2}dx+o\left(
1\right)  =\frac{1}{2H\left(  2H-1\right)  }+o\left(  1\right)  .
\]
We conclude that%
\begin{align*}
&  \mathbf{E}\left(  T_{2}^{2}\right)  +o\left(  N^{2H-2}\right) =\frac{cst.N^{2H-2}}{H(2H-1)} \\
& \times \sum_{q_{1},r_{1},q_{2},r_{2}=0}^{\ell
}b_{q_{1}}b_{r_{1}}b_{q_{2}}b_{r_{2}}\int_{[0,1]^{4}}dudvdu^{^{\prime}%
}dv^{^{\prime}}|u-v-q_{1}+r_{1}|^{2H^{^{\prime}}-2}|u^{^{\prime}}-v^{^{\prime
}}-q_{2}+r_{2}|^{2H^{^{\prime}}-2}.%
\end{align*}
Using the fact that
\begin{align*}
&  \int_{[0,1]^{2}}|u-v-q+r|^{2H^{^{\prime}}-2}dudv\\
&  =\frac{1}{2H^{^{\prime}}(2H^{^{\prime}}-1)}\left[  |1+q-r|^{2H^{^{\prime}}%
}+|1-q+r|^{2H^{^{\prime}}}-2|q-r|^{2H^{^{\prime}}}\right]
\end{align*}
the proposition follows.

\subsection{Proof of Proposition \ref{propET4}}

\subsubsection{The term $\mathbf{E}\left(  T_{4,(1)}^{2}\right)  $}

We have%
\begin{align*}
\mathbf{E}\left(  T_{4,(1)}^{2}\right)   &  =\frac{4N^{4H}}{c(H)^{2}%
(N-\ell)^{2}}4!\sum_{i,j=\ell}^{N-1}\langle C_{i}\otimes C_{i},C_{j}\otimes
C_{j}\rangle_{L^{2}([0,1]^{4})}\\
&  =\frac{4N^{4H}}{c(H)^{2}(N-\ell)^{2}}4!\sum_{i,j=\ell}^{N-1}\left\vert
\langle C_{i},C_{j}\rangle_{L^{2}([0,1]^{2})}\right\vert ^{2}%
\end{align*}
The scalar product computes as
\begin{align*}
\langle C_{i},C_{j}\rangle_{L^{2}([0,1]^{2})}  &  =\left\langle \sum
_{q=0}^{\ell}\alpha_{q}L_{\frac{i-q}{N}},\;\sum_{r=0}^{\ell}\alpha_{r}%
L_{\frac{j-r}{N}}\right\rangle _{L^{2}([0,1]^{2})}\\
&  =\int_{0}^{1}\int_{0}^{1}\left(  \sum_{q=0}^{\ell}\alpha_{q}L_{\frac
{i-q}{N}}(y_{1},y_{2})\right)  \left(  \sum_{r=0}^{\ell}\alpha_{r}%
L_{\frac{j-r}{N}}(y_{1},y_{2})\right)  dy_{1}dy_{2}\\
&  =d(H)^{2}\sum_{q,r=0}^{\ell}\alpha_{q}\alpha_{r}\int_{0}^{1}\int_{0}%
^{1}\left[  \int_{y_{1}\vee y_{2}}^{\frac{i-q}{N}}\frac{\partial K^{H^{\prime
}}}{\partial u}(u,y_{1})\frac{\partial K^{H^{\prime}}}{\partial u}%
(u,y_{2})du\right] \\
&  \;\;\times\left[  \int_{y_{1}\vee y_{2}}^{\frac{j-r}{N}}\frac{\partial
K^{H^{\prime}}}{\partial v}(v,y_{1})\frac{\partial K^{H^{\prime}}}{\partial
v}(v,y_{2})dv\right]  dy_{1}dy_{2}\\
&  =d(H)^{2}\sum_{q,r=0}^{\ell}\alpha_{q}\alpha_{r}\int_{0}^{\frac{i-q}{N}%
}\int_{0}^{\frac{j-r}{N}}\left(  \int_{0}^{u\wedge v}\frac{\partial
K^{H^{\prime}}}{\partial u}(u,y_{1})\frac{\partial K^{H^{\prime}}}{\partial
v}(v,y_{1})dy_{1}\right)  ^{2}dudv\\
&  =\alpha(H)^{2}\;d(H)^{2}\sum_{q,r=0}^{\ell}\alpha_{q}\alpha_{r}\int
_{0}^{\frac{i-q}{N}}\int_{0}^{\frac{j-r}{N}}|u-v|^{2H-2}dudv
\end{align*}
where $\alpha(H)=\frac{H(H+1)}{2}=H^{^{\prime}}(2H^{^{\prime}}-1)$ and
\begin{equation}
\int_{0}^{\frac{i-q}{N}}\int_{0}^{\frac{j-r}{N}}|u-v|^{2H-2}dudv=\frac
{1}{H(2H-1)}\left[  \left\vert \frac{i-q}{N}\right\vert ^{2H}+\left\vert
\frac{j-r}{N}\right\vert ^{2H}-\left\vert \frac{j-i+q-r}{N}\right\vert
^{2H}\right]  \label{int}%
\end{equation}
Using the fact that $\frac{\alpha(H)^{2}\;d(H)^{2}}{H(2H-1)}=\frac{1}{2}$ and
(\ref{int}) the scalar product becomes
\begin{align*}
\langle C_{i},C_{j}\rangle_{L^{2}([0,1]^{2})}  &  =\frac{\alpha(H)^{2}%
\;d(H)^{2}}{H(2H-1)}\;\sum_{q,r=0}^{\ell}\alpha_{q}\alpha_{r}\left[
\left\vert \frac{i-q}{N}\right\vert ^{2H}+\left\vert \frac{j-r}{N}\right\vert
^{2H}-\left\vert \frac{j-i+q-r}{N}\right\vert ^{2H}\right] \\
&  =\frac{1}{2}\sum_{q,r=0}^{\ell}\alpha_{q}\alpha_{r}\left[  \left\vert
\frac{i-q}{N}\right\vert ^{2H}+\left\vert \frac{j-r}{N}\right\vert
^{2H}-\left\vert \frac{j-i+q-r}{N}\right\vert ^{2H}\right] \\
&  =\frac{1}{2}\;\Bigg[\left(  \sum_{q=0}^{\ell}\alpha_{q}\left\vert
\frac{i-q}{N}\right\vert ^{2H}\right)  \left(  \sum_{r=0}^{\ell}\alpha
_{r}\right)  +\left(  \sum_{r=0}^{\ell}\alpha_{r}\left\vert \frac{j-r}%
{N}\right\vert ^{2H}\right)  \left(  \sum_{q=0}^{\ell}\alpha_{q}\right) \\
&  -\sum_{q,r=0}^{\ell}\alpha_{q}\alpha_{r}\left\vert \frac{i-j+q-r}%
{N}\right\vert ^{2H}\Bigg]\\
&  =-\frac{1}{2}\sum_{q,r=0}^{\ell}\alpha_{q}\alpha_{r}\left\vert
\frac{i-j+q-r}{N}\right\vert ^{2H}=\pi_{H}^{\alpha}(i-j).
\end{align*}
The last equality is true since $\sum_{q=0}^{\ell}\alpha_{q}=0$ by the filter
definition. Therefore, we have
\begin{align*}
&  \sum_{i,j=\ell}^{N-1}\left\vert \langle C_{i},C_{j}\rangle_{L^{2}%
([0,1]^{2})}\right\vert ^{2}\\
&  =\frac{1}{4}\;\sum_{i,j=\ell}^{N-1}\left\vert \sum_{q,r=0}^{\ell}\alpha
_{q}\alpha_{r}\left\vert \frac{i-j+q-r}{N}\right\vert ^{2H}\right\vert
^{2}=\frac{1}{4}\;\sum_{i=\ell}^{N-1}\sum_{k=0}^{N-2}\left\vert \sum
_{q,r=0}^{\ell}\alpha_{q}\alpha_{r}\left\vert \frac{k+q-r}{N}\right\vert
^{2H}\right\vert ^{2}\\
&  =\frac{N^{-4H}}{4}(N-\ell-1)\left\vert \sum_{q,r=0}^{\ell}\alpha_{q}%
\alpha_{r}|q-r|^{2H}\right\vert ^{2}+\frac{1}{4}\;\sum_{i=\ell}^{N-1}%
\sum_{k=1}^{N-2}\left\vert \sum_{q,r=0}^{\ell}\alpha_{q}\alpha_{r}\left\vert
\frac{k+q-r}{N}\right\vert ^{2H}\right\vert ^{2}\\
&  =c(H)^{2}\frac{N^{-4H}(N-\ell-1)}{4}+\frac{1}{4}\;\sum_{k=0}^{N-2}%
(N-k-2)\left\vert \sum_{q,r=0}^{\ell}\alpha_{q}\alpha_{r}\left\vert
\frac{k+q-r}{N}\right\vert ^{2H}\right\vert ^{2}\\
&  =c(H)^{2}\frac{(N-l-1)N^{-4H}}{4}+\frac{N^{-4H+1}}{4}\;\sum_{k=0}%
^{N-2}\left\vert \sum_{q,r=0}^{\ell}\alpha_{q}\alpha_{r}\left\vert
k+q-r\right\vert ^{2H}\right\vert ^{2}\\
&  -2\;\frac{N^{-4H}}{4}\;\sum_{k=0}^{N-2}\left\vert \sum_{q,r=0}^{\ell}%
\alpha_{q}\alpha_{r}\left\vert k+q-r\right\vert ^{2H}\right\vert ^{2}%
+\frac{N^{-4H}}{4}\;\sum_{k=0}^{N-2}k\left\vert \sum_{q,r=0}^{\ell}\alpha
_{q}\alpha_{r}\left\vert k+q-r\right\vert ^{2H}\right\vert ^{2}.
\end{align*}
At this point we need the next lemma to estimate the behavior of the above
quantity. This lemma is the key point which implies the fact that the longer
variation statistics has, in the case when the observed process is the
fractional Brownian motion, a Gaussian limit without any restriction on $H$
(see \cite{GuLe}).

\begin{itemize}
\item
\begin{lemma}
For all $H\in(0,1)$, we have that

\begin{description}
\item[(i)] $\sum_{k=1}^{\infty}\left\vert \sum_{q,r=0}^{\ell}\alpha_{q}%
\alpha_{r}|k+q-r|^{2H}\right\vert ^{2}<+\infty$

\item[(ii)] $\sum_{k=1}^{\infty}k\left\vert \sum_{q,r=0}^{\ell}\alpha
_{q}\alpha_{r}|k+q-r|^{2H}\right\vert ^{2}<+\infty .$
\end{description}
\end{lemma}

\begin{proof}
\emph{Proof of (i).} Let $f(x)=\sum_{q,r=0}^{\ell}\alpha_{q}\alpha_{r}\left(
1+(q-r)x\right)  ^{2H}$, so the summand can be written as
\[
\sum_{q,r=0}^{\ell}\alpha_{q}\alpha_{r}|k+q-r|^{2H}=k^{2H}f\left(  \frac{1}%
{k}\right).
\]
Using a Taylor expansion at $x_{0}=0$ for the function $f(x)$ we get that
\[
\left(  1+(q-r)x\right)  ^{2H}\approx1+2H(q-r)x+\ldots+\frac{2H(2H-1)\ldots
(2H-n+1)}{n!}(q-r)^{n}x^{n}.%
\]
For small $x$ we observe that the function $f(x)$ is asymptotically equivalent
to
\[
2H(2H-1)\ldots(2H-(p-1))x^{2p},
\]
where $p$ is the order of the filter. Therefore, the general term of the
series is equivalent to
\[
(2H)^{2}(2H-1)^{2}\ldots(2H-(p-1))^{2}k^{4H-4p}%
\]
Therefore for all $H<p-\frac{1}{4}$ the series converges to a constant
depending only on $H$. Due to our choice for the order of the filter $p\geq2$,
we obtain the desired result.

\emph{Proof of (ii).} Similarly as before, we can write the general term of
the series as
\begin{align*}
k\left\vert \sum_{q,r=0}^{\ell}\alpha_{q}\alpha_{r}|k+q-r|^{2H}\right\vert
^{2}  &  =k\left\vert k^{2H}f\left(  \frac{1}{k}\right)  \right\vert ^{2}\\
&  \approx(2H)^{2}(2H-1)^{2}\ldots(2H-(p-1))^{2}k^{4H-4p-1}%
\end{align*}
Therefore for all $H<p$ the series converges to a constant depending only on
$H$.
\end{proof}

\noindent Combining all the above we have
\begin{align*}
&  \mathbf{E}\left(  T_{4,(1)}^{2}\right)  =\frac{4\;N^{4H}}{c(H)^{2}%
(N-\ell)^{2}}4!\sum_{i,j=1}^{N}\left\vert \langle C_{i},C_{j}\rangle
_{L^{2}([0,1]^{2})}\right\vert ^{2}\\
&  =\frac{4\;N^{4H}}{c(H)^{2}(N-\ell)^{2}}4!\Bigg [\frac{1}{4}\;c(H)^{2}%
(N-\ell-1)N^{-4H}+\frac{N^{-4H+1}}{4}\;\sum_{k=0}^{N-2}\left\vert \sum
_{q,r=0}^{\ell}\alpha_{q}\alpha_{r}\left\vert k+q-r\right\vert ^{2H}%
\right\vert ^{2}\\
&  -2\;\frac{N^{-4H}}{4}\;\sum_{k=0}^{N-2}\left\vert \sum_{q,r=0}^{\ell}%
\alpha_{q}\alpha_{r}\left\vert k+q-r\right\vert ^{2H}\right\vert ^{2}%
+\frac{N^{-4H}}{4}\;\sum_{k=0}^{N-2}k\left\vert \sum_{q,r=0}^{\ell}\alpha
_{q}\alpha_{r}\left\vert k+q-r\right\vert ^{2H}\right\vert ^{2}\Bigg]\\
&  =\frac{4!}{c(H)^{2}}\Bigg [c(H)^{2}\;\frac{N-\ell-1}{(N-\ell)^{2}}+\left(
\frac{N^{1}}{(N-\ell)^{2}}-2\frac{1}{(N-\ell)^{2}}\right)  \;\sum_{k=0}%
^{N-2}\left\vert \sum_{q,r=0}^{\ell}\alpha_{q}\alpha_{r}\left\vert
k+q-r\right\vert ^{2H}\right\vert ^{2}\\
&  +\frac{1}{(N-\ell)^{2}}\;\sum_{k=0}^{N-2}k\left\vert \sum_{q,r=0}^{\ell
}\alpha_{q}\alpha_{r}\left\vert k+q-r\right\vert ^{2H}\right\vert ^{2}\Bigg]\\
&  =\frac{4!}{c(H)^{2}}\Bigg [c(H)^{2}\;\left(  \frac{N}{(N-\ell)^{2}}%
-\frac{l+1}{(N-\ell)^{2}}\right)  +\frac{N-2}{(N-\ell)^{2}}\;\sum_{k=0}%
^{N-2}\left\vert \sum_{q,r=0}^{\ell}\alpha_{q}\alpha_{r}\left\vert
k+q-r\right\vert ^{2H}\right\vert ^{2}\\
&  +\frac{1}{(N-\ell)^{2}}\;\sum_{k=0}^{N-2}k\left\vert \sum_{q,r=0}^{\ell
}\alpha_{q}\alpha_{r}\left\vert k+q-r\right\vert ^{2H}\right\vert ^{2}\Bigg]\\
&  \approx\frac{4!}{c(H)^{2}}\Bigg [c(H)^{2}\;\left(  N^{-1}-(\ell
+1)N^{-2}\right)  +\left(  N^{-1}-2N^{-2}\right)  \;\sum_{k=0}^{N-2}\left\vert
\sum_{q,r=0}^{\ell}\alpha_{q}\alpha_{r}\left\vert k+q-r\right\vert
^{2H}\right\vert ^{2}\\
&  +N^{-2}\;\sum_{k=0}^{N-2}k\left\vert \sum_{q,r=0}^{\ell}\alpha_{q}%
\alpha_{r}\left\vert k+q-r\right\vert ^{2H}\right\vert ^{2}\Bigg].
\end{align*}
Since the leading term is of order $N^{-1}$ we have that
\[
\mathbf{E}\left(  T_{4,(1)}^{2}\right)  \simeq4!\;c(H)^{-2}N^{-1}\left[
c(H)^{2}+\sum_{k=0}^{N-2}\left\vert \sum_{q,r=0}^{\ell}\alpha_{q}\alpha
_{r}\left\vert k+q-r\right\vert ^{2H}\right\vert ^{2}\right]  .
\]
If we define the correlation function of the filtered process as
\[
\rho_{H}^{\alpha}(k)=\frac{\pi_{H}^{\alpha}(k)}{\pi_{H}^{\alpha}(0)}%
=\frac{\sum_{q,r=0}^{\ell}\alpha_{q}\alpha_{r}\left\vert k+q-r\right\vert
^{2H}}{c(H)}%
\]
we can express the asymptotic variance $\lim_{N\rightarrow\infty}%
N~\mathbf{E}\left(  T_{4,(1)}^{2}\right)  $ in terms of a series involving
$\rho_{H}^{\alpha}(k)$.
\end{itemize}

\subsubsection{The term $\mathbf{E}\left(  T_{4,(2)}^{2}\right)  $}

In order to handle this term we use the alternate expression (\ref{SeeNose})
of $C_{i}$. Therefore, following similar calculations as in the $T_{2}$ case
we get that
\begin{align*}
\mathbf{E}\left(  T_{4,(2)}^{2}\right)   &  =\frac{c_{4,H}^{(1)}}{(N-\ell
)^{2}}\sum_{q_{1},q_{2},r_{1},r_{1}=0}^{\ell}b_{q_{1}}b_{q_{2}}b_{r_{1}%
}b_{r_{2}}\int_{[0,1]^{4}}dudvdu^{\prime}dv^{\prime}\\
&  \times\sum_{i,j=\ell}^{N-1}\bigg[\left\vert u-v+i-j-q_{1}+r_{1}\right\vert
^{2H^{\prime}-2}\left\vert u^{\prime}-v^{\prime}+i-j-q_{2}+r_{2}\right\vert
^{2H^{\prime}-2}\\
&  \left\vert u-u^{\prime}+i-j-q_{1}+q_{2}\right\vert ^{2H^{\prime}%
-2}\left\vert v-v^{\prime}+i-j-r_{1}+r_{2}\right\vert ^{2H^{\prime}-2}\bigg]\\
&  =\frac{c_{4,H}^{(2)}}{(N-\ell)^{2}}\sum_{q_{1},q_{2},r_{1},r_{1}=0}^{\ell
}b_{q_{1}}b_{q_{2}}b_{r_{1}}b_{r_{2}}\int_{[0,1]^{4}}dudvdu^{\prime}%
dv^{\prime}\\
&  \times\sum_{i=\ell}^{N-1}\sum_{k=0}^{N-\ell-i}\bigg[\left\vert
u-v+k-q_{1}+r_{1}\right\vert ^{2H^{\prime}-2}\left\vert u^{\prime}-v^{\prime
}+k-q_{2}+r_{2}\right\vert ^{2H^{\prime}-2}\\
&  \left\vert u-u^{\prime}+k-q_{1}+q_{2}\right\vert ^{2H^{\prime}-2}\left\vert
v-v^{\prime}+k-r_{1}+r_{2}\right\vert ^{2H^{\prime}-2}\bigg]\\
&  =\frac{c_{4,H}^{(3)}}{(N-\ell)^{2}}\sum_{q_{1},q_{2},r_{1},r_{1}=0}^{\ell
}b_{q_{1}}b_{q_{2}}b_{r_{1}}b_{r_{2}}\int_{[0,1]^{4}}dudvdu^{\prime}%
dv^{\prime}\\
&  \times\sum_{k=0}^{N-\ell}(N-k-1)\bigg[\left\vert u-v+k-q_{1}+r_{1}%
\right\vert ^{2H^{\prime}-2}\left\vert u^{\prime}-v^{\prime}+k-q_{2}%
+r_{2}\right\vert ^{2H^{\prime}-2}\\
&  \left\vert u-u^{\prime}+k-q_{1}+q_{2}\right\vert ^{2H^{\prime}-2}\left\vert
v-v^{\prime}+k-r_{1}+r_{2}\right\vert ^{2H^{\prime}-2}\bigg].
\end{align*}
We study the convergence of the above series as $N\rightarrow\infty$
\begin{align*}
&  \sum_{k=0}^{N-1}(N-k-1)\bigg[\left\vert u-v+k-q_{1}+r_{1}\right\vert
^{2H^{\prime}-2}\left\vert u^{\prime}-v^{\prime}+k-q_{2}+r_{2}\right\vert
^{2H^{\prime}-2}\\
&  \left\vert u-u^{\prime}+k-q_{1}+q_{2}\right\vert ^{2H^{\prime}-2}\left\vert
v-v^{\prime}+k-r_{1}+r_{2}\right\vert ^{2H^{\prime}-2}\bigg]\\
&  =(N-1)\sum_{k=0}^{N-1}\bigg[\left\vert u-v+k-q_{1}+r_{1}\right\vert
^{2H^{\prime}-2}\left\vert u^{\prime}-v^{\prime}+k-q_{2}+r_{2}\right\vert
^{2H^{\prime}-2}\\
&  \left\vert u-u^{\prime}+k-q_{1}+q_{2}\right\vert ^{2H^{\prime}-2}\left\vert
v-v^{\prime}+k-r_{1}+r_{2}\right\vert ^{2H^{\prime}-2}\bigg]\\
&  -\sum_{k=0}^{N-1}k\bigg[\left\vert u-v+k-q_{1}+r_{1}\right\vert
^{2H^{\prime}-2}\left\vert u^{\prime}-v^{\prime}+k-q_{2}+r_{2}\right\vert
^{2H^{\prime}-2}\\
&  \left\vert u-u^{\prime}+k-q_{1}+q_{2}\right\vert ^{2H^{\prime}-2}\left\vert
v-v^{\prime}+k-r_{1}+r_{2}\right\vert ^{2H^{\prime}-2}\bigg]\\
&  :=(I)\;+\;(II).
\end{align*}
Therefore the general term of the series is asymptotically equivalent to
\begin{align*}
&  \left(  \frac{(2H^{\prime}-2)\ldots(2H^{\prime}-2p-1)}{(2p)!}\right)
^{4}(u-v-q_{1}+r_{1})^{2p}\;(u^{\prime}-v^{\prime}-q_{2}+r_{2})^{2p}\\
&  \hspace*{1.5in}\cdot(u-u^{\prime}-q_{1}+q_{2})^{2p}\;(v-v^{\prime}%
-r_{1}+r_{2})^{2p}\;\;k^{4H-4-8p},
\end{align*}
which converges for all $H\in(\frac{1}{2},1)$. We treat the second series (II)
in the same way and we get that it is asymptotically equivalent to
$cst.\;k^{4H-4-8p}$. Combining all the above we have
\begin{align*}
\mathbf{E}\left(  T_{4,(2)}^{2}\right)   &  =\frac{c_{4,H}^{^{\prime}}%
}{(N-\ell)^{2}}\sum_{q_{1},q_{2},r_{1},r_{1}=0}^{\ell}b_{q_{1}}b_{q_{2}%
}b_{r_{1}}b_{r_{2}}\int_{[0,1]^{4}}dudvdu^{\prime}dv^{\prime}\\
&  \bigg\{(N-\ell)\sum_{k=\ell}^{N-1}\bigg[\left\vert u-v+k-q_{1}%
+r_{1}\right\vert ^{2H^{\prime}-2}\left\vert u^{\prime}-v^{\prime}%
+k-q_{2}+r_{2}\right\vert ^{2H^{\prime}-2}\\
&  \left\vert u-u^{\prime}+k-q_{1}+q_{2}\right\vert ^{2H^{\prime}-2}\left\vert
v-v^{\prime}+k-r_{1}+r_{2}\right\vert ^{2H^{\prime}-2}\bigg]\\
&  -\sum_{k=\ell}^{N-1}k\bigg[\left\vert u-v+k-q_{1}+r_{1}\right\vert
^{2H^{\prime}-2}\left\vert u^{\prime}-v^{\prime}+k-q_{2}+r_{2}\right\vert
^{2H^{\prime}-2}\\
&  \left\vert u-u^{\prime}+k-q_{1}+q_{2}\right\vert ^{2H^{\prime}-2}\left\vert
v-v^{\prime}+k-r_{1}+r_{2}\right\vert ^{2H^{\prime}-2}\bigg]\bigg \}.
\end{align*}
The leading term in $\mathbf{E}\left(  T_{4,(2)}^{2}\right)  $ is of order
$N^{-1}$ and the constant computes as
\begin{align*}
\tau_{1,H}  &  =\sum_{k=\ell}^{\infty}\sum_{q_{1},q_{2},r_{1},r_{1}=0}^{\ell
}b_{q_{1}}b_{q_{2}}b_{r_{1}}b_{r_{2}}\int_{[0,1]^{4}}dudvdu^{\prime}%
dv^{\prime}\\
&  \bigg[\left\vert u-v+k-q_{1}+r_{1}\right\vert ^{2H^{\prime}-2}\left\vert
u^{\prime}-v^{\prime}+k-q_{2}+r_{2}\right\vert ^{2H^{\prime}-2}\\
&  \left\vert u-u^{\prime}+k-q_{1}+q_{2}\right\vert ^{2H^{\prime}-2}\left\vert
v-v^{\prime}+k-r_{1}+r_{2}\right\vert ^{2H^{\prime}-2}\bigg].
\end{align*}

Therefore, combining the two terms we get the statement of the proposition.

\subsection{End of proof of Theorem \ref{thmT4}}

Recall that we only need to show that for $\tau=1,2,3$ the terms $||g_{N}\otimes_{\tau} g_{N}||^{2}_{L^{2}([0,1]^{8-2\tau})}$ converge to 0 as $N$ tends to infinity.

\begin{itemize}

\item \emph{Term} for $\tau=1$.

\begin{eqnarray*}
J_{1} &=& \left( \frac{4 N^{4H+1}}{c_{1,H}c(H)^{2}(N-\ell)^{2}} \right)^{2}
\sum_{i,j,m,n=\ell}^{N-1}      \left \langle (C_{i}\otimes C_{i}) \otimes_{1} (C_{j}\otimes C_{j}) ,
(C_{m}\otimes C_{m}) \otimes_{1} (C_{n}\otimes C_{n}) \right \rangle \\
&=& \left( \frac{4 N^{4H+1}}{c_{1,H}c(H)^{2}(N-\ell)^{2}} \right)^{2}
\sum_{i,j,m,n=\ell}^{N-1}      \langle C_{i}, C_{m}\rangle_{L^{2}([0,1]^{2})}\;\langle C_{j}, C_{n}\rangle_{L^{2}([0,1]^{2})}\\
&& \times \langle C_{i}\otimes_{1} C_{j}, C_{m} \otimes_{1} C_{n} \rangle_{L^{2}([0,1]^{2})}.
\end{eqnarray*}

Thus, we have
\begin{align*}
&  J_{1}  =\\
&  \leq cst.\frac{N^{8H+2}}{(N-\ell)^{4}}\frac{1}{N^{4}}\sum_{i,j,m,n=\ell
}^{N-1}\;\;\sum_{q_{1},r_{1},q_{2},r_{2},q_{3},r_{3},q_{4},r_{4}=0}^{\ell
}b_{q_{1}}b_{r_{1}}b_{q_{2}}b_{r_{2}}b_{q_{3}}b_{r_{3}}b_{q_{4}}b_{r_{4}}\\
&  \times\left\vert \frac{i-m+q_{1}-r_{1}}{N}\right\vert ^{2H}\;\left\vert
\frac{j-n+q_{2}-r_{2}}{N}\right\vert ^{2H}\\
&  \times\bigg[\int_{[0,1]^{4}}dudvdu^{\prime}dv^{\prime}\left\vert
\frac{u-v+i-j-q_{3}+r_{3}}{N}\right\vert ^{2H^{\prime}-2}\left\vert
\frac{u^{\prime}-v^{\prime}+m-n-q_{4}+r_{4}}{N}\right\vert ^{2H^{\prime}-2}\\
&  \times\left\vert \frac{u-u^{\prime}+i-m-q_{3}+q_{4}}{N}\right\vert
^{2H^{\prime}-2}\left\vert \frac{v-v^{\prime}+j-n+r_{3}+r_{4}}{N}\right\vert
^{2H^{\prime}-2}\bigg]\\
&  \leq cst.\frac{N^{2}}{(N-\ell)^{4}}\sum_{i,j,m,n=\ell}^{N-1}\;\;\sum
_{q_{1},r_{1},q_{2},r_{2},q_{3},r_{3},q_{4},r_{4}=0}^{\ell}b_{q_{1}}b_{r_{1}%
}b_{q_{2}}b_{r_{2}}b_{q_{3}}b_{r_{3}}b_{q_{4}}b_{r_{4}}\\
&  \times\left\vert i-m+q_{1}-r_{1}\right\vert ^{2H}\;\left\vert
j-n+q_{2}-r_{2}\right\vert ^{2H}\\
&  \times\bigg[\int_{[0,1]^{4}}dudvdu^{\prime}dv^{\prime}|u-v+i-j-q_{3}%
+r_{3}|^{2H^{\prime}-2}|u^{\prime}-v^{\prime}+m-n-q_{4}+r_{4}|^{2H^{\prime}%
-2}\\
&  \times|u-u^{\prime}+i-m-q_{3}+q_{4}|^{2H^{\prime}-2}|v-v^{\prime}%
+j-n+r_{3}+r_{4}|^{2H^{\prime}-2}\bigg].
\end{align*}
As in the computations for $T_{4,(2)}$ we can show that the above series
converges and thus $  J_{1}  ={\mathcal{O}}%
(N^{-2})$, which implies that for all $H\in(\frac{1}{2},1)$
\[
\lim_{N\rightarrow\infty}  J_{1}  =0.
\]

\item \emph{Term} for $\tau=2$

\begin{eqnarray*}
J_{2} &=& \left( \frac{4 N^{4H+1}}{c_{1,H}c(H)^{2}(N-\ell)^{2}} \right)^{2}
\sum_{i,j,m,n=\ell}^{N-1}      \left \langle (C_{i}\otimes C_{i}) \otimes_{2} (C_{j}\otimes C_{j}) ,
(C_{m}\otimes C_{m}) \otimes_{2} (C_{n}\otimes C_{n}) \right \rangle \\
&=& \left( \frac{4 N^{4H+1}}{c_{1,H}c(H)^{2}(N-\ell)^{2}} \right)^{2}
\sum_{i,j,m,n=\ell}^{N-1}      \langle C_{i}, C_{j}\rangle_{L^{2}([0,1]^{2})}\;\langle C_{m}, C_{n}\rangle_{L^{2}([0,1]^{2})}\\
&&  \times\langle C_{i}, C_{m}\rangle_{L^{2}([0,1]^{2})}\;\langle C_{j}, C_{n}\rangle_{L^{2}([0,1]^{2})}.
\end{eqnarray*}

\end{itemize}%

\begin{align*}
J_{2} &  \leq cst.\frac{N^{8H+2}%
}{(N-\ell)^{4}}\sum_{i,j,m,n=\ell}^{N-1}\langle C_{i},C_{j}\rangle
_{L^{2}[0,1]^{2}}\langle C_{i},C_{m}\rangle_{L^{2}[0,1]^{2}}\langle
C_{m},C_{n}\rangle_{L^{2}[0,1]^{2}}\langle C_{j},C_{n}\rangle_{L^{2}[0,1]^{2}%
}\\
&  =cst.\frac{N^{8H+2}}{(N-\ell)^{4}}\sum_{i,j,m,n=\ell}^{N-1}\sum_{q_{1}%
q_{2}q_{3}q_{4}=0}^{\ell}\alpha_{q_{1}}\alpha_{q_{2}}\alpha_{q_{3}}%
\alpha_{q_{4}}\left\vert \frac{i-j+q_{1}-q_{2}}{N}\right\vert ^{2H}\\
&  \times\left\vert \frac{i-m+q_{1}-q_{3}}{N}\right\vert ^{2H}\left\vert
\frac{m-n+q_{3}-q_{4}}{N}\right\vert ^{2H}\left\vert \frac{j-n+q_{2}-q_{4}}%
{N}\right\vert ^{2H}\\
&  =cst.\frac{N^{2}}{(N-\ell)^{4}}\sum_{i,j,m,n=\ell}^{N-1}\sum_{q_{1}%
q_{2}q_{3}q_{4}=0}^{\ell}\alpha_{q_{1}}\alpha_{q_{2}}\alpha_{q_{3}}%
\alpha_{q_{4}}\left\vert i-j+q_{1}-q_{2}\right\vert ^{2H}\\
&  \times\left\vert i-m+q_{1}-q_{3}\right\vert ^{2H}\left\vert m-n+q_{3}%
-q_{4}\right\vert ^{2H}\left\vert j-n+q_{2}-q_{4}\right\vert ^{2H}.
\end{align*}
The series converges for all $H\in(1/2,1)$, so the whole term is of order
${\mathcal{O}}(N^{-2})$ which means that goes to zero as $N\rightarrow\infty$.

\begin{itemize}
\item \emph{Term} for $\tau=3$.
\begin{eqnarray*}
J_{3} &=& \left( \frac{4 N^{4H+1}}{c_{1,H}c(H)^{2}(N-\ell)^{2}} \right)^{2}
\sum_{i,j,m,n=\ell}^{N-1}      \left \langle (C_{i}\otimes C_{i}) \otimes_{3} (C_{j}\otimes C_{j}) ,
(C_{m}\otimes C_{m}) \otimes_{3} (C_{n}\otimes C_{n}) \right \rangle \\
&=& \left( \frac{4 N^{4H+1}}{c_{1,H}c(H)^{2}(N-\ell)^{2}} \right)^{2}
\sum_{i,j,m,n=\ell}^{N-1}      \langle C_{i}, C_{j}\rangle_{L^{2}([0,1]^{2})}\;\langle C_{m}, C_{n}\rangle_{L^{2}([0,1]^{2})}\\
&& \times \langle C_{i}\otimes_{1} C_{j}, C_{m} \otimes_{1} C_{n} \rangle.
\end{eqnarray*}

With similar computations as in the case of $T_{4}$
we conclude that $J_{3} ={\mathcal{O}}(N^{-2})$.
\end{itemize}

\subsection{Proof of Theorem \ref{thmV}}

According to our previous computations we can write
\begin{align*}
&  f_{N}(y_{1},y_{2})=\\
&  =\frac{8N^{2H}}{c(H)(N-\ell)}\sum_{i=\ell}^{N-1}(C_{i}\otimes_{1}%
C_{i})(y_{1},y_{2})\\
&  =\frac{8d(H)^{2}\alpha(H)}{c(H)}\frac{N^{2H}}{(N-\ell)}\sum_{i=\ell}%
^{N-1}\sum_{q,r=0}^{\ell}b_{q}b_{r}1_{[0,\frac{i-q+1}{N}]}(y_{1}%
)1_{[0,\frac{i-r+1}{N}]}(y_{2})\\
&  \times\int_{I_{i_{q}}}\int_{I_{i_{r}}}dudv|u-v|^{2H^{\prime}-2}%
\;\partial_{1}K^{H^{\prime}}(u,y_{1})\partial_{1}K^{H^{\prime}}(v,y_{2})
\end{align*}
Let us show first that we can reduce this function to the interval $y_{1}%
\in\lbrack0,\frac{i-q}{N}]$ and $y_{2}\in\lbrack0,\frac{i-r}{N}]$. We will
show that if $y_{1}\in I_{i_{q}},y_{2}\in\lbrack0,\frac{i-r}{N}]$ (and
similarly for the situations $y_{1}\in\lbrack0,\frac{i-q}{N}],y_{2}\in
I_{i_{r}}$ and $y_{1}\in I_{i_{q}},y_{2}\in I_{i_{r}}$) the corresponding
terms goes to zero as $N\rightarrow\infty$. We have, due to the fact that the
intervals $I_{i_{q}}$ are disjoint,
\begin{align*}
&  \Vert\frac{N^{1-H}N^{2H}}{(N-\ell)}\sum_{i=\ell}^{N-1}\sum_{q,r=0}^{\ell
}b_{q}b_{r}1_{I_{i_{q}}}(y_{1})1_{[0,\frac{i-r}{N}]}(y_{2}) \\
& \int_{I_{i_{q}}
}\int_{I_{i_{r}}}dudv|u-v|^{2H^{\prime}-2}\;\partial_{1}K^{H^{\prime}}%
(u,y_{1})\partial_{1}K^{H^{\prime}}(v,y_{2})\Vert_{L^{2}([0,1]^{2})}^{2}\\
&  =\frac{N^{2+2H}}{(N-\ell)^{2}}\sum_{i=\ell}^{N}\sum_{q_{1},r_{1}%
,q_{2},r_{2}=0}^{\ell}b_{q_{1}}b_{r_{1}}b_{q_{2}}b_{r_{2}}\int_{I_{i_{q_{1}}}%
}\int_{I_{i_{r_{1}}}}\int_{I_{i_{q_{2}}}}\int_{I_{i_{r_{2}}}}dv^{\prime
}du^{\prime}dvdu\\
&  \times\left(  |u-v|\cdot|u^{\prime}-v^{\prime}|\cdot|u-u^{\prime}%
|\cdot|v-v^{\prime}|\right)  ^{2H^{\prime}-2}\\
&  =\frac{N^{2+2H}}{(N-\ell)^{2}}\frac{1}{N^{4}}\frac{1}{N^{4(2H^{\prime}-2)}%
}\sum_{i=\ell}^{N}\sum_{q_{1},r_{1},q_{2},r_{2}=0}^{\ell}b_{q_{1}}b_{r_{1}%
}b_{q_{2}}b_{r_{2}}\int_{[0,1]^{4}}dudvdu^{^{\prime}}dv^{^{\prime}}\\
&  \times|u-v-q_{1}+r_{1}|^{2H^{^{\prime}}-2}|u^{^{\prime}}-v^{^{\prime}%
}-q_{2}+r_{2}|^{2H^{^{\prime}}-2}\\
& |u-u^{^{\prime}}-q_{1}+q_{2}|^{2H^{^{\prime}%
}-2}|v-v^{^{\prime}}-r_{1}+r_{2}|^{2H^{^{\prime}}-2}\asymp{N^{1-2H}}%
\end{align*}
which tends to zero because $2H>1$.

This proves the following asymptotic equivalence in $L^{2}([0,1]^{2})$:%
\begin{align*}
f_{N}(y_{1},y_{2})  &  \simeq\frac{8d(H)^{2}\alpha(H)}{c(H)}\frac{N^{2H}%
}{(N-\ell)}\sum_{i=\ell}^{N-1}\sum_{q,r=0}^{\ell}b_{q}b_{r}1_{[0,\frac{i-q}%
{N}]}(y_{1})1_{[0,\frac{i-r}{N}]}(y_{2})\\
&  \times\int_{I_{i_{q}}}\int_{I_{i_{r}}}dudv|u-v|^{2H^{\prime}-2}%
\;\partial_{1}K^{H^{\prime}}(u,y_{1})\partial_{1}K^{H^{\prime}}(v,y_{2}).
\end{align*}
We will show that the above term, normalize by $\frac{N^{1-H}}{\sqrt{c_{2,H}}%
}$, converges pointwise for $y_{1},\;y_{2}\in\lbrack0,1]$ to the kernel of the
Rosenblatt random variable.

On the interval $I_{i_{q}}\times I_{i_{r}}$ we may attemp to replace the
evaluation of $\partial_{1}K^{H^{\prime}}$ at $u$ and $v$ by setting
$u=(i-q)/N$ and $v=(i-r)/N$. More precisely, we can write%
\begin{align*}
&  \partial_{1}K^{H^{\prime}}(u,y_{1})\partial_{1}K^{H^{\prime}}%
(v,y_{2})=\left(  \partial_{1}K^{H^{\prime}}(u,y_{1})-\partial_{1}%
K^{H^{\prime}}(\frac{i-q}{N},y_{1})\right)  \partial_{1}K^{H^{\prime}}%
(v,y_{2})\\
&  +\partial_{1}K^{H^{\prime}}(\frac{i-q}{N},y_{1})\left(  \partial
_{1}K^{H^{\prime}}(v,y_{2})-\partial_{1}K^{H^{\prime}}-\partial_{1}%
K^{H^{\prime}}(\frac{i-r}{N},y_{2})\right)
\end{align*}
and all the above summand above can be treated in the same manner. For the
first one, using the definition of the derivative of $K^{H^{\prime}}$ with
respect to the first variable, we get for any $y_{1}\in\lbrack0,\frac{i-q}%
{N}]$,
\begin{align*}
&  \partial_{1}K^{H^{\prime}}(u,y_{1})-\partial_{1}K^{H^{\prime}}(\frac
{i-q}{N},y_{1})\\
&  =c_{H}y_{1}^{\frac{1}{2}-H}\left(  (u-y_{1})^{H-\frac{3}{2}}u^{H-\frac
{1}{2}}-\left(  \frac{i-q}{N}-y_{1}\right)  ^{H-\frac{3}{2}}(\frac{i-q}%
{N})^{H-\frac{1}{2}}\right) \\
&  \leq c_{H}y_{1}^{\frac{1}{2}-H}\left(  \frac{i-q}{N}-y_{1}\right)
^{H-\frac{3}{2}}\left(  u^{H-\frac{1}{2}}-(\frac{i-q}{N})^{H-\frac{1}{2}%
}\right) \\
&  \leq c_{H}y_{1}^{\frac{1}{2}-H}\left(  \frac{i-q}{N}-y_{1}\right)
^{H-\frac{3}{2}}(u-(\frac{i-q}{N}))^{H-\frac{1}{2}}\\
&  \leq c_{H}N^{\frac{1}{2}-H}y_{1}^{\frac{1}{2}-H}\left(  \frac{i-q}{N}%
-y_{1}\right)  ^{H-\frac{3}{2}}%
\end{align*}
and for any $y_{2}\in\lbrack0,\frac{i-r}{N}]$
\begin{align*}
\partial_{1}K^{H^{\prime}}(v,y_{2})  &  =c_{H}y_{2}^{\frac{1}{2}-H}%
(v-y_{2})^{H-\frac{3}{2}}v^{H-\frac{1}{2}}\\
&  \leq c_{H}y_{2}^{\frac{1}{2}-H}\left(  \frac{i-r}{N}-y_{1}\right)
^{H-\frac{3}{2}}(\frac{i-r+1}{N})^{H-\frac{1}{2}}.
\end{align*}
As a consequence of the above estimates,
\begin{align}
&  N^{1-H}\frac{N^{2H}}{N-\ell}\sum_{i=\ell}^{N-1}\sum_{q,r=0}^{\ell}%
b_{q}b_{r}1_{[0,\frac{i-q}{N}]}(y_{1})1_{[0,\frac{i-r}{N}]}(y_{2})\nonumber\\
&  \times\int_{I_{i_{q}}}\int_{I_{i_{q}}}dvdu|u-v|^{2H^{\prime}-2}\left(
\partial_{1}K^{H^{\prime}}(u,y_{1})-\partial_{1}K^{H^{\prime}}(\frac{i-q}%
{N},y_{1})\right)  \partial_{1}K^{H^{\prime}}(v,y_{2})\nonumber\\
&  \leq cN^{\frac{1}{2}-H}\frac{N^{1+H}}{N-\ell}\sum_{i=\ell}^{N-1}%
\sum_{q,r=0}^{\ell}b_{q}b_{r}1_{[0,\frac{i-q}{N}]}(y_{1})1_{[0,\frac{i-r}{N}%
]}(y_{2})\nonumber\\
&  \times\left(  \frac{i-q}{N}-y_{1}\right)  ^{H-\frac{3}{2}}\left(
\frac{i-r}{N}-y_{2}\right)  ^{H-\frac{3}{2}}(\frac{i-r+1}{N})^{H-\frac{1}{2}%
}\int_{I_{i_{q}}}\int_{I_{i_{q}}}dvdu|u-v|^{2H^{\prime}-2}\nonumber\\
&  \leq cN^{\frac{1}{2}-H}\frac{1}{N-\ell}\sum_{i=\ell}^{N-1}\sum
_{q,r=0}^{\ell}b_{q}b_{r}1_{[0,\frac{i-q}{N}]}(y_{1})1_{[0,\frac{i-r}{N}%
]}(y_{2})\nonumber\\
&  \times\left(  \frac{i-q}{N}-y_{1}\right)  ^{H-\frac{3}{2}}\left(
\frac{i-r}{N}-y_{2}\right)  ^{H-\frac{3}{2}}(\frac{i-r+1}{N})^{H-\frac{1}{2}}.
\label{1dif}%
\end{align}

\noindent The quantity $\frac{1}{N-\ell}\sum_{i=\ell}^{N-1}1_{[0,\frac{i-q}{N}]}%
(y_{1})1_{[0,\frac{i-r}{N}]}(y_{2})\left(  \frac{i-q}{N}-y_{1}\right)
^{H-\frac{3}{2}}\left(  \frac{i-r}{N}-y_{2}\right)  ^{H-\frac{3}{2}}%
(\frac{i-r+1}{N})^{H-\frac{1}{2}}$ is comparable, for large $N,$ to the
integral $\int_{y_{1}\vee y_{2}}^{1}(u-y_{1})^{H-\frac{3}{2}}(u-y_{2}%
)^{H-\frac{3}{2}}u^{H-\frac{1}{2}}$ and the term $N^{\frac{1}{2}-H}$ in front
gives the convergence to zero of (\ref{1dif}) for any fixed $y_{1},y_{2}$.

This means we have proved the following pointwise asymptotically equivalent
for $f_{N}(y_{1},y_{2})$:%
\begin{align*}
\frac{N^{1-H}}{\sqrt{c_{2,H}}}f_{N}(y_{1},y_{2})  &  \simeq\frac
{8d(H)^{2}\alpha(H)}{\sqrt{c_{2,H}}\;c(H)}\frac{N^{1+H}}{(N-\ell)}\sum
_{i=\ell}^{N-1}\sum_{q,r=0}^{\ell}1_{[0,\frac{i-q}{N}]}(y_{1})1_{[0,\frac
{i-r}{N}]}(y_{2})b_{q}b_{r}\\
&  \times\;\partial_{1}K^{H^{\prime}}(\frac{i-q}{N},y_{1})\partial
_{1}K^{H^{\prime}}(\frac{i-r}{N},y_{2})\;\int_{I_{i_{q}}}\int_{I_{i_{r}}%
}dudv|u-v|^{2H^{\prime}-2}.%
\end{align*}
Recall that
\[
\int_{I_{i_{q}}}\int_{I_{i_{r}}}dvdu|u-v|^{2H^{\prime}-2}=\frac{N^{-(1+H)}%
}{2H^{\prime}(2H^{\prime}-1)}\left\{  |1-q+r|^{2H^{\prime}}%
+|1+q-r|^{2H^{\prime}}-2|q-r|^{2H^{\prime}}\right\}  .
\]
Thus we get
\begin{align*}
&  \frac{N^{1-H}}{\sqrt{c_{2,H}}}f_{N}(y_{1},y_{2})\\
&  \simeq\frac{8d(H)^{2}\alpha(H)}{c_{2,H}\;c(H)}\sum_{q,r=0}^{\ell}b_{q}%
b_{r}\;\left\{  |1-q+r|^{2H^{\prime}}+|1+q-r|^{2H^{\prime}}-2|q-r|^{2H^{\prime
}}\right\} \\
&  \times\frac{1}{(N-\ell)}\sum_{i=\ell}^{N-1}\partial_{1}K^{H^{\prime}}%
(\frac{i-q}{N},y_{1})\partial_{1}K^{H^{\prime}}(\frac{i-r}{N},y_{2}).
\end{align*}
Further, we can ignore the terms $q/N$ and $r/N$ in comparison with $i/N$ in
the last line above, and thus invoke a Riemann sum approximation, which proves
that, for every $y_{1},y_{2}\in(0,1)^{2}$
\begin{align*}
&  \lim_{N\rightarrow\infty}\frac{N^{1-H}}{c_{2,H}}f_{N}(y_{1},y_{2})\\
&  =\frac{8d(H)^{2}\alpha(H)}{c_{2,H}\;c(H)}\sum_{q,r=0}^{\ell}b_{q}%
b_{r}\;\left\{  |1-q+r|^{2H^{\prime}}+|1+q-r|^{2H^{\prime}}-2|q-r|^{2H^{\prime
}}\right\} \\
&  \frac{1}{(N-\ell)}\lim_{N\rightarrow\infty}\sum_{i=\ell}^{N-1}\partial
_{1}K^{H^{\prime}}(\frac{i-q}{N},y_{1})\partial_{1}K^{H^{\prime}}(\frac
{i-r}{N},y_{2})\\
&  =d(H)\int_{y_{1}\vee y_{2}}\partial_{1}K^{H^{\prime}}(\frac{u}{N}%
,y_{1})\partial_{1}K^{H^{\prime}}(\frac{u}{N},y_{2})du\\
&  =L_{1}(y_{1},y_{2}).
\end{align*}
To finish the proof it suffices to check that $N^{1-H}f_{N}$ is a Cauchy
sequence in $L^{2}([0,1]^{2})$. Up to a constant depending on $H$ we have that
for all $M$, $N$,%
\begin{align*}
&  ||N^{1-H}f_{N}-M^{1-H}f_{M}||_{L^{2}([0,1]^{2})}^{2}\\
&  =N^{2-2H}||f_{N}||_{L^{2}([0,1]^{2})}^{2}+M^{2-2H}||f_{M}||_{L^{2}%
([0,1]^{2})}^{2}-2N^{1-H}M^{1-H}\langle f_{N},f_{M}\rangle_{L^{2}([0,1]^{2}%
)}\\
&  =cst.\frac{N^{2H+2}}{(N-\ell)^{2}}\sum_{i,j=\ell}^{N-1}\sum_{q_{1}%
,r_{1},q_{2},r_{2}=0}^{\ell}b_{q_{1}}b_{r_{1}}b_{q_{2}}b_{r_{2}}%
\int_{I_{i_{q_{1}}}^{N}}\int_{I_{i_{r_{1}}}^{N}}\int_{I_{j_{q_{2}}}^{N}}%
\int_{I_{j_{r_{2}}}^{N}}dudvdu^{\prime}dv^{\prime}\\
&  \times|u-v|^{2H^{\prime}-2}|u^{\prime}-v^{\prime}|^{2H^{\prime}%
-2}|u-u^{\prime}|^{2H^{\prime}-2}|v-v^{\prime}|^{2H^{\prime}-2}\\
&  +cst.\frac{M^{2H+2}}{(M-\ell)^{2}}\sum_{i,j=\ell}^{M-1}\sum_{q_{1}%
,r_{1},q_{2},r_{2}=0}^{\ell}b_{q_{1}}b_{r_{1}}b_{q_{2}}b_{r_{2}}%
\int_{I_{i_{q_{1}}}^{M}}\int_{I_{i_{r_{1}}}^{M}}\int_{I_{j_{q_{2}}}^{M}}%
\int_{I_{j_{r_{2}}}^{M}}dudvdu^{\prime}dv^{\prime}\\
&  \times|u-v|^{2H^{\prime}-2}|u^{\prime}-v^{\prime}|^{2H^{\prime}%
-2}|u-u^{\prime}|^{2H^{\prime}-2}|v-v^{\prime}|^{2H^{\prime}-2}\\
&  -cst.\frac{M^{1+H}N^{1+H}}{(M-\ell)(N-\ell)}\sum_{i=\ell}^{N-1}\sum
_{j=\ell}^{M-1}\sum_{q_{1},r_{1},q_{2},r_{2}=0}^{\ell}b_{q_{1}}b_{r_{1}%
}b_{q_{2}}b_{r_{2}}\int_{I_{i_{q_{1}}}^{N}}\int_{I_{i_{r_{1}}}^{N}}%
\int_{I_{j_{q_{2}}}^{M}}\int_{I_{j_{r_{2}}}^{M}}dudvdu^{\prime}dv^{\prime}\\
&  \times|u-v|^{2H^{\prime}-2}|u^{\prime}-v^{\prime}|^{2H^{\prime}%
-2}|u-u^{\prime}|^{2H^{\prime}-2}|v-v^{\prime}|^{2H^{\prime}-2}.%
\end{align*}
The first two terms have already been studied and will converge to the same
constant as $M,N\rightarrow\infty$. Concerning the inner product, by making
the usual change of variable we have
\begin{align*}
&  \frac{(MN)^{H+1}}{(M-\ell)(N-\ell)}\frac{(NM)^{2H^{\prime}-2}}{N^{2}M^{2}%
}\sum_{i=\ell}^{N-1}\sum_{j=\ell}^{M-1}\sum_{q_{1},r_{1},q_{2},r_{2}=0}^{\ell
}\int_{[0,1]^{4}}dudvdu^{\prime}dv^{\prime}\\
&  \times|u-v-q_{1}+r_{1}|^{2H^{\prime}-2}|u^{\prime}-v^{\prime}-q_{3}%
+r_{3}|^{2H^{\prime}-2}\\
&  \times\left\vert \frac{u}{N}-\frac{u^{\prime}}{M}+\frac{i}{N}-\frac{j}%
{N}-\frac{q_{1}}{N}+\frac{q_{2}}{N}\right\vert ^{2H^{\prime}-2}\left\vert
\frac{v}{N}-\frac{v^{\prime}}{M}+\frac{i}{N}-\frac{j}{N}-\frac{r_{1}}{N}%
+\frac{r_{2}}{N}\right\vert ^{2H^{\prime}-2}.%
\end{align*}
For large $i,j$ we can ignore the terms $\frac{u}{N}$, $\frac{u^{\prime}}{N}$,
$\frac{q_{1}}{N}$, etc., compared to $\frac{i}{N}$ and $\frac{j}{N}$.
Therefore, the above quantity is a Riemann sum that converges to the same
constant as the squared terms, as $M,N\rightarrow\infty$. This finishes the
proof of the theorem.

\subsection{Proof of Theorem \ref{thmstatconv}}

We wish to show that, as $N\rightarrow\infty$,%
\[
E:=\mathbf{E}\left[  \left(  Z\left(  1\right)  -2c_{2,H}^{-1/2}N^{1-H}\left(
\hat{H}_{N}-H\right)  \log N\right)  ^{2}\right]  \rightarrow0.
\]
A minor technical difficulty occurs when $V_{N}$ is not small. We deal with
this as follows. We decompose the above expectation $E$ according to whether
or not $\left\vert V_{N}\right\vert \leq1/2$: we have $E=E_{1}+E_{2}$ where%
\[
E_{1}=\mathbf{E}\left[  \mathbf{1}_{\left\vert V_{N}\right\vert >1/2}\left(
Z\left(  1\right)  -2c_{2,H}^{-1/2}N^{1-H}\left(  \hat{H}_{N}-H\right)  \log
N\right)  ^{2}\right]  .
\]
Dealing with this term first, Schwarz's and Minkowski's inequalities yields%
\[
E_{1}\leq2\mathbf{P}^{1/2}\left[  \left\vert V_{N}\right\vert >1/2\right]
\left(  \mathbf{E}^{1/2}\left[  Z\left(  1\right)  ^{4}\right]  +4c_{2,H}%
^{-1}N^{2-2H}\log^{2}N~\mathbf{E}^{1/2}\left[  \left(  \hat{H}_{N}-H\right)
^{4}\right]  \right)  .
\]
Since $\hat{H}_{N}$ is bounded, the sum of the two rooted expectation terms
above is bounded above by a constant multiple of $N^{2-2H}$. Therefore to deal
with $E_{1}$, one only needs to show that $\mathbf{P}\left[  \left\vert
V_{N}\right\vert >1/2\right]  \ll N^{-4+4H}$. It is well known that any random
variable $X$ which can be written as a finite sum of Wiener chaos terms up to
order $q$ satisfies, for any integer $n$, $\mathbf{E}\left[  X^{2n}\right]
\leq K_{n,q}\left(  \mathbf{E}\left[  X^{2}\right]  \right)  ^{n}$ where
$K_{n,q}$ depends only on $n$ and $q$. This can be proved iteratively by using
formula (\ref{prod}), for instance. Therefore, since $V_{N}$ is a sum of terms
in the second and $4$th chaos ($q=4$), by Chebyshev's inequality, and using
Theorem \ref{th}, with $N$ large enough,%
\begin{align*}
\mathbf{P}\left[  \left\vert V_{N}\right\vert >1/2\right]   &  \leq
4^{n}\mathbf{E}\left[  \left\vert V_{N}\right\vert ^{2n}\right]  \leq
4^{n}c_{n,4}\left(  \mathbf{E}\left[  \left\vert V_{N}\right\vert ^{2}\right]
\right)  ^{n}\\
&  \leq8^{n}K_{n,4}c_{2,H}^{n}N^{2Hn-2n}.
\end{align*}
It is thus sufficient to choose $n=3$ to guarantee that $E_{1}\rightarrow0$.

We now only need to study $E_{2}$. We invoke the mean value theorem to express
$\left(  \hat{H}_{N}-H\right)$  $\log N$ more explicitly. For any $x,y\in
\lbrack1/2,1]$, there exists $\xi\in(x,y)$ such that
\[
\log\frac{c\left(  x\right)  }{c\left(  y\right)  }=\left(  x-y\right)
\left(  \log c\right)  ^{\prime}\left(  \xi\right)  .
\]
Here the function $\left(  \log c\right)  ^{\prime}$ is bounded on $[1/2,1]$,
because $c^{\prime}$ is bounded and $c$ is bounded below. Therefore, denoting
by $\xi_{N}\in\lbrack1/2,1]$ the value corresponding to $x=H$ and $y=\hat
{H}_{N}$, and using line (\ref{thmconsistentline}) in the proof of Theorem
\ref{thmconsistent}, we can write%
\[
\log\left(  1+V_{N}\right)  =\left(  \hat{H}_{N}-H\right)  \left(  2\log
N+\left(  \log c\right)  ^{\prime}\left(  \xi_{N}\right)  \right)
\]
and thus%
\[
\left(  \hat{H}_{N}-H\right)  \left(  2\log N\right)  =\log\left(
1+V_{N}\right)  -\frac{\log\left(  1+V_{N}\right)  }{2\log N+\left(  \log
c\right)  ^{\prime}\left(  \xi_{N}\right)  }.
\]
Since $\left\vert \left(  \log c\right)  ^{\prime}\left(  \xi_{N}\right)
\right\vert $ is bounded (by a non-random value), by choosing $N$ large
enough, an upper bound for the last fraction above, in absolute value, is
$2V_{N}/\log N$. Therefore (using Minkowski's inequality),%
\begin{align}
\sqrt{E_{2}}  &  =\mathbf{E}^{1/2}\left[  \mathbf{1}_{\left\vert
V_{N}\right\vert \leq1/2}\left(  Z\left(  1\right)  -2c_{2,H}^{-1/2}%
N^{1-H}\left(  \hat{H}_{N}-H\right)  \log N\right)  ^{2}\right] \nonumber\\
&  \leq\mathbf{E}^{1/2}\left[  \mathbf{1}_{\left\vert V_{N}\right\vert
\leq1/2}\left(  Z\left(  1\right)  -c_{2,H}^{-1/2}N^{1-H}\left(  \log\left(
1+V_{N}\right)  \right)  \right)  ^{2}\right] \label{line2}\\
&  +\mathbf{E}^{1/2}\left[  \mathbf{1}_{\left\vert V_{N}\right\vert \leq
1/2}\left(  2c_{2,H}^{-1/2}N^{1-H}V_{N}/\log N\right)  ^{2}\right]  .
\label{line3}%
\end{align}

By Theorem \ref{th}, the term in line (\ref{line3}) is bounded above by
$1/\log^{2}N$, and thus converges to $0$. For the term in line (\ref{line2}),
because of the indicator $\mathbf{1}_{\left\vert V_{N}\right\vert \leq1/2}$,
we use the fact that when $\left\vert x\right\vert \leq1/2$, we have
$\left\vert x-\log\left(  1+x\right)  \right\vert \leq x^{2}$. Thus this line
is bounded above by
\begin{align}
&  \mathbf{E}^{1/2}\left[  \mathbf{1}_{\left\vert V_{N}\right\vert \leq
1/2}\left(  Z\left(  1\right)  -c_{2,H}^{-1/2}N^{1-H}V_{N}\right)  ^{2}\right]
\label{line4}\\
&  +\mathbf{E}^{1/2}\left[  \mathbf{1}_{\left\vert V_{N}\right\vert \leq
1/2}\left(  c_{2,H}^{-1/2}N^{1-H}\left\vert V_{N}\right\vert ^{2}\right)
^{2}\right] . \label{line5}%
\end{align}
The term in line (\ref{line4}) converges to $0$ by Theorem \ref{thmV}.
Finally, by Theorem \ref{th} again, and the earlier statement about higher
powers of random variables with finite chaos expansions, the term in line
(\ref{line5}) is of order $N^{2H-2}$, and therefore converges to $0$ as well.
This proves that $E_{2}$ converges to $0$, finishing the proof of the theorem.

\subsection{Proof of Theorem \ref{ThmStat}}

It is sufficient to prove that
\[
\lim_{N\rightarrow\infty}\mathbf{E}\left[  \left\vert \left(  N^{1-\hat{H}%
_{N}}-N^{1-H}\right)  \left(  \hat{H}_{N}-H\right)  \log N\right\vert \right]
=0.
\]
We decompose the probability space depending on whether $\hat{H}_{N}$ is far
or not from its mean. For a fixed value $\varepsilon>0$ it is convenient to
define the event
\[
D=\left\{  \hat{H}_{N}>\varepsilon+2H-1\right\}.
\]
We have
\begin{align*}
&  \mathbf{E}\left[  \left\vert \left(  N^{1-\hat{H}_{N}}-N^{1-H}\right)
\left(  \hat{H}_{N}-H\right)  \log N\right\vert \right]  =\\
&  =\mathbf{E}\left[  \mathbf{1}_{D}\left\vert \left(  N^{1-\hat{H}_{N}%
}-N^{1-H}\right)  \left(  \hat{H}_{N}-H\right)  \log N\right\vert \right]  +\\
&  +\mathbf{E}\left[  \mathbf{1}_{D^{c}}\left\vert \left(  N^{1-\hat{H}_{N}%
}-N^{1-H}\right)  \left(  \hat{H}_{N}-H\right)  \log N\right\vert \right] \\
&  :=A+B.
\end{align*}

\begin{proof}
\begin{description}
\item[Term A] :\newline Introduce the notation $x=\max\left(  1-H,1-\hat
{H}_{N}\right)  $ and $y=\min\left(  1-H,1-\hat{H}_{N}\right)  $.
\begin{align*}
\left\vert N^{1-\hat{H}_{N}}-N^{1-H}\right\vert  &  =e^{x\log N}-e^{y\log
N}=e^{y\log N}\left(  e^{(x-y)\log N}-1\right) \\
&  \leq N^{y}(\log N)(x-y)N^{x-y}=2\log NN^{x}\left\vert H-\hat{H}%
_{N}\right\vert \\
&  =\log NN^{x}\left\vert H-\hat{H}_{N}\right\vert.
\end{align*}
Thus,
\begin{align*}
A  &  \leq\mathbf{E}\left[  \mathbf{1}_{D}N^{x}\left\vert H-\hat{H}%
_{N}\right\vert ^{2}\log^{2}N\right] \\
&  =\mathbf{E}\left[  N^{x-(2-2H)}\mathbf{1}_{D}N^{2-2H}\left\vert H-\hat
{H}_{N}\right\vert ^{2}\log^{2}N\right]
\end{align*}
Now, choose $\varepsilon\in(0,1-H)$. In this case, if $\omega\in D$ and
$x=1-H$, we get $x-(2-2H)=-x<-\varepsilon$. On the other hand, for $\omega\in
D$ and $x=1-\hat{H}_{N}$ we get $x-(2-2H)=2-2\hat{H}_{N}-(2-2H)<-\varepsilon$.
In conclusion, on $D$, $x-(2-2H)<-\varepsilon$ which implies immediately
\[
A\leq N^{-\varepsilon}\mathbf{E}\left[  N^{2-2H}\left\vert \hat{H}%
_{N}-H\right\vert ^{2}\log^{2}N\right].
\]
and since the last expectation is bounded
\[
\lim_{N\rightarrow\infty}A=0.
\]

\item[Term B] :\newline Now, let $\omega\in D^{c}$ then $H-\hat{H}%
_{N}>1-H-\varepsilon$. Since $\varepsilon<1-H$ it implies $H>\hat{H}_{N}$.
Consequently, it is not sufficient to bound $\left\vert N^{1-\hat{H}_{N}%
}-N^{1-H}\right\vert $ above by $N^{1-\hat{H}_{N}}$. In the same fashion we
bound $\left\vert \hat{H}-H\right\vert $ above by $H$. Using H\"{o}lder's
inequality with powers $\frac{1}{4}$ and $\frac{3}{4}$
\begin{align*}
B  &  \leq H\log N\mathbf{E}\left[  \mathbf{1}_{D^{c}}N^{1-\hat{H}_{N}}\right]
\\
&  \leq H\log N\left[  P(D^{c})\right]  ^{3/4}\left(  \mathbf{E}\left[
N^{(1-\hat{H}_{N})4}\right]  \right)  ^{1/4}.%
\end{align*}

By Chebyshev's inequality, we have
\begin{equation}
\mathbf{P}^{3/4}\left[  D^{c}\right]  \leq\frac{\mathbf{E}^{3/4}\left[
\left\vert H-\hat{H}\right\vert ^{2}\right]  }{\left(  1-H-\varepsilon\right)
^{3/2}}\leq cN^{-3\left(  2-2H\right)  /4}%
\end{equation}
for some constant $c$ depending only $H$. Dealing with the other term in the
upper bound for $B$ is a little less obvious. We must return to the definition
of $\hat{H}$. We have
\[
1+V_{N}=N^{2\left(  H-\hat{H}\right)  }=N^{4\left(  H-\hat{H}\right)
}=N^{4\left(  1-\hat{H}\right)  }N^{-4\left(  1-H\right)  }.
\]
Therefore,
\[
\mathbf{E}^{1/4}\left[  N^{\left(  1-\hat{H}\right)  4}\right]  \leq
N^{1-H}\mathbf{E}^{1/\left(  4\right)  }\left[  1+V_{N}\right]  \leq2N^{1-H}.
\]
Finally, we get
\[
B\leq2Hc\left(  \log N\right)  N^{-\left(  1-H\right)  }.
\]
Finally, $B$ goes to $0$ as $N\rightarrow\infty$. This finishes the proof of
the theorem.
\end{description}
\end{proof}

\end{document}